\documentclass[12pt]{amsart}
\usepackage{epic,eepic}
\usepackage{graphicx}
\usepackage{amsthm}
\usepackage{amsfonts}
\usepackage{amssymb}
\usepackage{latexsym}
\usepackage{amscd}

\def\phi{\varphi}
\def\1#1{^{(#1)}}
\newcommand{\I}[1]{I\1{#1}}
\newcommand{\area}[1]{\operatorname{area}(#1)}

\DeclareMathOperator{\diam}{\mathrm{diam}}

\def\Nu{\mathcal{V}}

\catcode`\@=12

\def\Empty{}
\newcommand\oplabel[1]{
  \def\OpArg{#1} \ifx \OpArg\Empty {} \else
  	\label{#1}
  \fi}
		
\long\def\realfig#1#2#3#4{
\begin{figure}[htbp]
\centerline{\includegraphics[width=#4]{#2}}
\caption[#1]{#3}
\oplabel{#1}
\end{figure}}

\newcommand{\comm}[1]{}

\input{psfig}
\newtheorem{thm}{Theorem}[section]
\newtheorem{cor}[thm]{Corollary}
\newtheorem{lem}[thm]{Lemma}
\newtheorem{prop}[thm]{Proposition}

\newenvironment{pf}{\proof[\proofname]}{\endproof}
\newenvironment{pf*}[1]{\proof[#1]}{\endproof}
\usepackage{euscript}

\usepackage[OT2,OT1]{fontenc}

\newcommand{\cal}[1]{{\mathcal #1}}

\theoremstyle{definition}
\newtheorem{defn}{Definition}[section]

\theoremstyle{remark}
\newtheorem{rem}{Remark}[section]

\newcommand{\dist}{\operatorname{dist}}

\newcommand{\cl}{\operatorname{cl}}
\renewcommand{\mod}{\operatorname{mod}}
\newcommand{\tl}{\tilde}
\newcommand{\wtl}{\widetilde}
\newcommand{\eps}{\epsilon}
\newcommand{\EE}{{\Ccal{E}}}

\newcommand{\ceq}{\,\displaystyle{\Large\mbox{$\sim$}}_{\text{\hspace{-14pt}\tiny conf}}\,}

\newcommand{\Cbb}[1]{{{\Bbb{#1}}}}
\newcommand{\Ccal}[1]{{{\cal{#1}}}}
\newcommand{\aaa}[1]{{{\mathbf{#1}}}}
\newcommand{\crit}{{{\aaa C}}}
\newcommand{\cu}{{{\aaa C}_U}}
\newcommand{\cur}{{{\aaa C}_U^\RR}}
\newcommand{\cv}{{{\aaa C}_V}}

\newcommand{\hol}{{\aaa H}}
\newcommand{\mfld}{{\aaa M}}
\newcommand{\bran}{{{\aaa X}}}

\newcommand{\abs}[1]{|{#1}|}

\newcommand{\sg}{\sigma}
\newcommand{\dl}{\delta}

\renewcommand{\Im}{\operatorname{Im}}

\numberwithin{equation}{section}
\newcommand{\thmref}[1]{Theorem~\ref{#1}}
\newcommand{\propref}[1]{Proposition~\ref{#1}}
\newcommand{\secref}[1]{\S\ref{#1}}
\newcommand{\lemref}[1]{Lemma~\ref{#1}}
\newcommand{\corref}[1]{Corollary~\ref{#1}}
\newcommand{\figref}[1]{Figure~\ref{#1}}

\newcommand{\cW}{{\cal W}}
\newcommand{\cM}{{\cal M}}

\newcommand{\cF}{{\cal F}}
\newcommand{\cG}{{\cal G}}

\newcommand{\cI}{{\cal I}}
\newcommand{\cP}{{\cal P}}
\newcommand{\cC}{{\cal C}}
\newcommand{\cH}{{\cal H}}
\newcommand{\cR}{{\cal R}}
\newcommand{\cL}{{\cal L}}
\newcommand{\cD}{{\cal D}}
\newcommand{\cE}{{\cal E}}

\newcommand{\CC}{{\Bbb C}}
\newcommand{\RR}{{\Bbb R}}
\newcommand{\TT}{{\Bbb T}}
\newcommand{\ZZ}{{\Bbb Z}}
\newcommand{\NN}{{\Bbb N}}
\newcommand{\DD}{{\Bbb D}}
\newcommand{\HH}{{\Bbb H}}
\newcommand{\QQ}{{\Bbb Q}}

\newcommand{\cren}{\cR_{\text cyl}}

\newcommand{\sm}{\setminus}

\begin{document}
\addtolength{\evensidemargin}{-0.7in}
\addtolength{\oddsidemargin}{-0.7in}

\title[Rigidity Problem]{The rigidity problem for analytic critical circle maps}
\author{\fbox{D. Khmelev} and M. Yampolsky}
\thanks{This paper was completed during the authors' stay at the Institut Henri Poincar{\'e}
during Fall 2003, and presented at the IHP at the same time.
\\
The second author was partially supported by an NSERC Discovery grant.}
\date{September 8, 2004}
\begin{abstract}
It is shown that if $f$ and $g$ are any two analytic critical circle mappings with the same irrational rotation number, then
the conjugacy that maps the critical point of $f$ to that of $g$ has regularity $C^{1+\alpha}$ at the critical point, with a
universal value of $\alpha>0$. As a consequence, a new proof of the hyperbolicity of the full renormalization horseshoe
of critical circle maps is given.
\end{abstract}
\maketitle


\noindent
{\bf Foreword.} In the last several decades, since the works of 
Mostow, Margulis, Sullivan, and others,  {\it rigidity problems}
occupy a central place in the theory of holomorphic dynamical systems. 
This type of problems is classical in dynamics: 
a rigidity theorem postulates that in a certain class of dynamical
systems equivalence (combinatorial, continuous, smooth, etc.) 
automatically has a higher regularity.
The dynamical systems considered in this paper are {\it critical circle
maps}, that is smooth homeomorphisms of the circle with a single
critical point having a cubic type. These maps have been a subject
of intensive study since the early 1980's as one of the two main examples
of universality in transition to chaos. In 1984 Yoccoz \cite{Yoc}
showed that any two such maps with the same irrational rotation number
are conjugate by a continuous change of coordinates
(this is a
generalization of the classical result of Denjoy for $C^2$ diffeomorphisms).
The central result of this paper is the following rigidity theorem:

\medskip
\noindent
{\bf Main Theorem}. {\sl For any two analytic critical circle maps with the same irrational rotation number,
the conjugacy, which maps the critical point of one to the critical point of the
other is $C^{1+\alpha}$-smooth at the critical point. Here $\alpha>0$ is a universal
constant}.

\medskip
\noindent
The theorem verifies the so-called
$C^{1+\alpha}$-rigidity conjecture for critical circle maps. The conjecture
has a history dating back to at least the early 1980's, as it is related to
the universality in critical circle maps, and together with the similar conjecture
for the unimodal maps of the interval it appears in the works of Feigenbaum, Lanford, Sullivan
and others. 
There was some recent progress in the study of the conjecture, as de~Faria and de~Melo \cite{dFdM2} used the
methods of McMullen \cite{McM2} to establish it for the
rotation numbers satisfying the Diophantine condition with the exponent $2$, usually referred to in the
subject as the bounded type condition. The understanding of the unbounded type requires
an analysis of critical circle mappings with an almost parabolic dynamics, which we have
carried out in this paper.

\section{Preliminaries}
\label{section preliminaries}

\noindent
{\bf Some notations.}
We use $\dist$ and $\diam$ to denote the Euclidean distance and diameter in $\Bbb C$, and we let
 $\area{D}$ denote the area of a domain $D\subset\CC$. 
We shall say that two real numbers $A$ and $B$ are {\it $K$-commensurable} for $K>1$ if
$K^{-1}|A|\leq |B|\leq K|A|$. We will use the notation $A\underset{K}{\sim}B$ 
in this case, omitting $K$ and writing simply $A\sim B$ when the commensurability
factor is universal.
The notation $D_r(z)$ will stand for the Euclidean disk with the center at $z\in\Bbb C$ and
radius $r$. The unit disk $D_1(0)$ will be denoted $\Bbb D$. 
The plane $({\Bbb C}\setminus{\Bbb R})\cup J$
with the parts of the real axis not contained in the interval $J\subset \Bbb R$
removed  will be denoted ${\Bbb C}_J$.
By the circle $\TT$ we understand the affine manifold $\RR/\ZZ$, it is naturally identified
with the unit circle $S^1=\partial\DD$. The real translation $x\mapsto x+\theta$
projects to the {\it rigid rotation by angle $\theta$} of the circle, $R_\theta:\TT\to\TT$.
For two points $a$ and $b$ in the circle $\TT$
which are not diametrically opposite, $[a,b]$ will denote the shorter
of the two arcs connecting them. As usual, $|[a,b]|$ will denote the length of the 
arc. For two points $a,b\in \Bbb R$, $[a,b]$ will denote the closed interval with
endpoints $a$, $b$ without specifying their order.
{\it The cylinder} in this paper, unless otherwise specified, will mean the affine manifold $\CC/\ZZ$.
Its {\it equator} is the circle $\{\Im z=0\}/\ZZ\subset \CC/\ZZ$.
A topological  annulus $A\subset \CC/\ZZ$ will be called an {\it equatorial annulus}, or an 
{\it equatorial neighborhood}, if it has a smooth boundary and contains the equator.

By ``smooth'' in this paper we will mean ``of class $C^\infty$'', unless another degree
of smoothness is specified. The notation``$C^\omega$'' will 
stand for ``real-analytic''.

For a piecewise-analytic dynamical system $\cF=\sqcup f_i$ in the Riemann sphere
we will call the dynamics generated by the iterates of $f_i$ and the univalent inverse branches
$f_i^{-1}$ the {\it complete dynamics of} $\cF$.

\subsection{Statements of the results.}
A critical circle map is an orientation preserving automorphism of
 $\TT$ of class $C^3$ with a single critical
point $c$. A further assumption is made that the critical point is of cubic type.
This means that for a lift $\bar f:{\Bbb R}\to {\Bbb R}$ of a critical circle  map $f$
with critical points at the integer translates of $\bar c$,
$$\bar f(x)-\bar f(\bar c)=(x-\bar c)^3(\operatorname{const}+O(x-\bar c)).$$
We note that all the renormalization results will hold true if in the above definition
``3'' as the order of smoothness and the order of the critical point is replaced
by any other odd number. To fix our ideas, we will always place the critical
point of $f$ at $0\in\TT$.

Being a homeomorphism of the circle, a critical circle map $f$  has a well-defined
rotation number, denoted $\rho(f)$. 
It is useful to represent $\rho(f)$ as a contined fraction with positive terms
\begin{equation}
\label{rotation-number}
\rho(f)=\cfrac{1}{r_0+\cfrac{1}{r_1+\cfrac{1}{r_2+\dotsb}}}
\end{equation}
Further on we will abbreviate this expression as $[r_0,r_1,r_2,\ldots]$
for typographical convenience. 
Note that the numbers $r_i$ are determined uniquely if and only if $\rho(f)$ is
irrational. In this case we shall say that $\rho(f)$ (or $f$ itself) is of the type bounded by $B$
if $\sup r_i\leq B$.

The main result of this paper is the following Theorem:

\begin{thm}
\label{main theorem}
There exists a universal constant $\alpha>0$ such that the following holds. Let $f_1$ and $f_2$
be two analytic critical circle maps with the same irrational rotation number. Denote
$\psi:\TT\to\TT$ the conjugacy $\psi\circ f_1\circ \psi^{-1}=f_2$ fixing the origin.
Then $\psi$ is $C^{1+\alpha}$ at the origin.

\end{thm}

\noindent
The above theorem should be seen as a generalization of the result of de~Faria and de~Melo
\cite{dFdM1,dFdM2}, who showed that when $\rho(f_i)$ is of a type bounded by some 
constant $B$, the conjugacy $\phi$ is globally $C^{1+\alpha}$ smooth, with $\alpha=\alpha(B)$.
An immediate corrolary is:

\begin{cor}
\label{ren conv}
The uniform distance between the successive renormalizations $\cR^n f_i$ decreases at a 
universal geometric rate.
\end{cor}

\noindent
We use this to obtain a new  proof of the main result of \cite{Ya4}:

\begin{thm}
\label{renormalization horseshoe}
The global renormalization horseshoe of the cylinder renormalization operator $\cren$
is uniformly hyperbolic, with one-dimensional unstable direction.
\end{thm}

\section{The geometry of the closest returns and renormalization}
\subsection{The dynamical partition of a critical circle map.}
Recall that an iterate $f^k(0)$ is called a {\it closest return} of the critical point if
the arc $[0,f^k(0)]$ contains no iterates $f^i(0)$ with $i<k$.
By a classical result of Poincar{\'e} every circle homeomorphism $f$ with
an irrational rotation number is semi-conjugate to the rigid rotation $R_{\rho(f)}$.
Moreover, Yoccoz \cite{Yoc} has shown that in the case when $f$ is a critical circle map,
the semi-conjugacy becomes a conjugacy, thus extending the classical result of Denjoy to
this case.
Poincar{\'e}'s result implies that 
the order of the points in an orbit of $f$ with an irrational rotation number 
is the same as that in an orbit of $R_{\rho(f)}$.
It follows, in particular, that if we denote $\{p_m/q_m\}$ the sequence of best rational 
approximations of $\rho(f)$ obtained as the truncated continued fractions
$p_m/q_m=[r_0,r_1,\ldots,r_{m-1}]$, then the iterates $\{f^{q_m}(0)\}$ 
are closest returns of $0$.
Set $I_m\equiv[0,f^{q_m}(0)]$. We will denote $I_m^i=f^i(I_m)$. An important combinatorial 
fact is that the collection of intervals
\begin{equation}
\label{dynamical partition}
\cP_m=\{I_m,I_m^1,\ldots,I_m^{q_{m+1}-1}\}\cup\{I_{m+1},I_{m+1}^1,\ldots,I_{m+1}^{q_m-1}\}
\end{equation}
covers the circle. Moreover, the intervals in the collection $\cP_m$ may overlap
only at the endpoints. Hence we will refer to $\cP_m$ as the {\it $m$-th dynamical
partition of $f$}. 
The geometry of the partition (\ref{dynamical partition}) is essential to our study.
Of crucial importance is the following {\it real {\em a priori} bound} of \'Swia\c\negthinspace tek and Herman:

\begin{thm}
\label{real bounds} There exists a universal constant $K>1$ such that the following holds.
Let $f:\TT\to\TT$ be a critical circle map with an irrational rotation
number. Then there exists $m_0=m_0(f)$ such that for all $m\geq m_0$ and every pair $I$, $J$
of adjacent atoms of the partition (\ref{dynamical partition}),
$$K^{-1}|J|\leq |I|\leq K|J|.$$
\end{thm}

\noindent
Utilizing the above bounds Herman has shown in \cite{H} (see \cite{dFdM1} for a published account):
\begin{thm}
Any two critical circle maps with the same irrational rotation numbers are quasisymmetrically conjugate.
\end{thm}

\subsection{Definition of renormalization of critical circle maps.}
An analogy with the universality phenomena in  statistical physics and with the 
already discovered Feigenbaum-Collett-Tresser universality in unimodal maps,
led the authors of \cite{FKS} and \cite{ORSS} to explain the existence of the universal
constants by introducing a renormalization operator acting on critical circle maps.
The definition is by no means straightforward. A detailed discussion may be found in
\cite{Ya3}. We need a supporting definition:

\begin{defn}
A  {\it  commuting pair} $\zeta=(\eta,\xi)$ consists of two 
$C^3$-smooth  orientation preserving interval homeomorphisms 
$\eta:I_\eta\to \eta(I_\eta),\;
\xi:I_{\xi}\to \xi(I_\xi)$, where
\begin{itemize}
\item[(I)]{$I_\eta=[0,\xi(0)],\; I_\xi=[\eta(0),0]$;}
\item[(II)]{Both $\eta$ and $\xi$ have homeomorphic extensions to interval
neighborhoods of their respective domains with the same degree of
smoothness, which commute, 
$\eta\circ\xi=\xi\circ\eta$;}
\item[(III)]{$\xi\circ\eta(0)\in I_\eta$;}
\item[(IV)]{$\eta'(x)\ne 0\ne \xi'(y) $, for all $x\in I_\eta\setminus\{0\}$,
 and all $y\in I_\xi\setminus\{0\}$.}
\end{itemize}
\end{defn}

\realfig{compair}{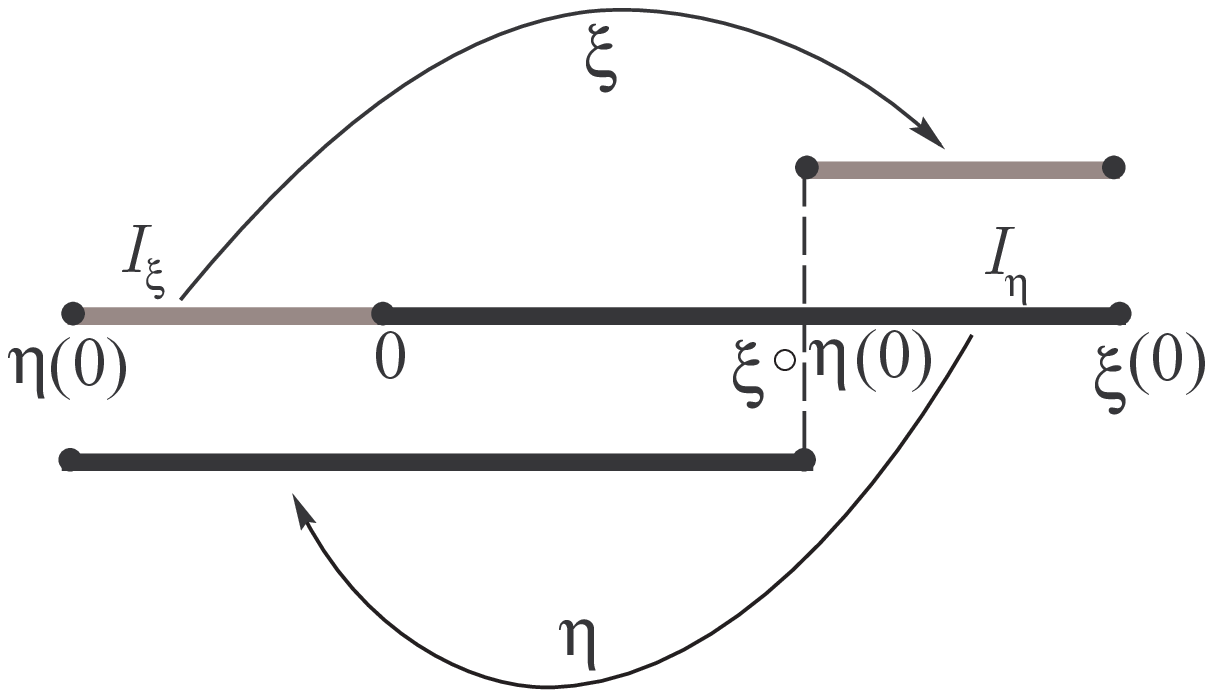}{A commuting pair}{8cm}

\noindent
The commutation condition allows one to iterate the extensions of the maps of a commuting pair.
It also allows us to perform the following glueing construction.
Given a critical commuting pair $\zeta=(\eta,\xi)$
we can regard the interval $I=[\eta(0),\xi\circ\eta(0)]$ as a circle, identifying 
$\eta(0)$ and $\xi\circ\eta(0)$ and define $f_\zeta:I\to I$ by 
$$f_\zeta=\left\{\begin{array}{l}
                    \eta\circ\xi(x) \text{ for }x\in [\eta(0),0]\\
                    \eta(x)\text{ for } x\in [0,\xi\circ\eta(0)]
                \end{array}\right.
$$
The mapping $\xi$ extends to a $C^3$-diffeomorphism of open
neighborhoods of  $\eta(0)$ and $\xi\circ\eta(0)$. Using it as a local chart we turn 
the interval $I$ into a closed one-dimensional manifold $M$. Condition (II) above implies that
the mapping $f_\zeta$ projects to a well-defined $C^3$-smooth homeomorphism $F_\zeta:M\to M$.
Identifying $M$ with the circle by a diffeomorphism $\phi:M\to \TT$
we recover a critical circle mapping $f^\phi=\phi\circ F_\zeta\circ\phi^{-1}$.
The critical circle mappings corresponding to two different choices
of $\phi$ are conjugated by a diffeomorphism, and thus we recovered
a $C^3$-smooth conjugacy class of circle mappings from a critical commuting
pair.

We can metrize the space of $C^r$-smooth commuting pairs considered modulo an affine conjugacy
as follows (see \cite{dFdM1}).
Let $\zeta_1=(\eta_1,\xi_1)$, $\zeta_2=(\eta_2,\xi_2)$
 be two such pairs, and denote $w_i:\CC\to\CC$ a M{\"o}bius transformation
which maps the ordered triple of points $\eta_i(0)$, $0$, $\xi_i(0)$ to $0$, $1/2$, $1$.
The {\it $C^r$-distance} between $\zeta_1$ and $\zeta_2$ is set to be 
$$\dist_{C^r}(\zeta_1,\zeta_2)=\max\{|\xi_1(0)/\eta_1(0)-\xi_2(0)/\eta_2(0)|,\dist_{C^r}(w_1\circ\zeta_1\circ w_1^{-1},w_2\circ\zeta_2\circ w_2^{-1})\}.$$

Let $f$ be a critical circle mapping, whose rotation number $\rho$
has a continued fraction expansion (\ref{rotation-number}) with
at least $m+1$ terms, and let $p_m/q_m=[r_0,\ldots,r_{m-1}]$. The pair of iterates $f^{q_{m+1}}$
and $f^{q_m}$ restricted to the circle arcs $I_m$ and $I_{m+1}$
correspondingly can be viewed as a critical commuting pair  in
the following way.
Let $\bar f$ be the lift of $f$ to the real line satisfying $\bar f '(0)=0$,
and $0<\bar f (0)<1$. For each $m>0$ let $\bar I_m\subset \Bbb R$ 
denote the closed 
interval adjacent to zero which projects down to the interval $I_m$.
Let $\tau :\Bbb R\to \Bbb R$ denote the translation $x\mapsto x+1$.
Let $\eta :\bar I_m\to \Bbb R$, $\xi:\bar I_{m+1}\to \Bbb R$ be given by
$\eta\equiv \tau^{-p_{m+1}}\circ\bar f^{q_{m+1}}$,
$\xi\equiv \tau^{-p_m}\circ\bar f^{q_m}$. Then the pair of maps
$(\eta|\bar I_m,\xi|\bar I_{m+1})$ forms a critical commuting pair
corresponding to $(f^{q_{m+1}}|I_m,f^{q_m}|I_{m+1})$.
Henceforth we shall  simply denote this commuting pair by
\begin{equation}
\label{real1}
(f^{q_{m+1}}|I_m,f^{q_m}|I_{m+1}).
\end{equation}
This allows us to readily identify the dynamics of
the above commuting pair with that of the underlying circle map,
at the cost of a minor abuse of notation.

Following \cite{dFdM1}, we say that the {\it height} $\chi(\zeta)$
of a critical commuting pair $\zeta=(\eta,\xi)$ is equal to $r$,
if 
$$0\in [\eta^r(\xi(0)),\eta^{r+1}(\xi(0))].$$
 If no such $r$ exists,
we set $\chi(\zeta)=\infty$, in this case the map $\eta|I_\eta$ has a 
fixed point.  For a pair $\zeta$ with $\chi(\zeta)=r<\infty$ one verifies directly that the
mappings $\eta|[0,\eta^r(\xi(0))]$ and $\eta^r\circ\xi|I_\xi$
again form a commuting pair.
For a commuting pair $\zeta=(\eta,\xi)$ we will denote by 
$\wtl\zeta$ the pair $(\wtl\eta|\wtl{I_\eta},\wtl\xi|\wtl{I_\xi})$
where tilde  means rescaling by the linear factor $\lambda=-{1\over |I_\eta|}$.

\begin{defn}
The {\it renormalization} of a real commuting pair $\zeta=(\eta,
\xi)$ is the commuting pair
\begin{center}
${\cal{R}}\zeta=(
\widetilde{\eta^r\circ\xi}|
 \widetilde{I_{\xi}},\; \widetilde\eta|\widetilde{[0,\eta^r(\xi(0))]}).$
\end{center}
\end{defn}

\noindent
The non-rescaled pair $(\eta^r\circ\xi|I_\xi,\eta|[0,\eta^r(\xi(0))])$ will be referred to as the 
{\it pre-renormalization} $p{\cal R}\zeta$ of the commuting pair $\zeta=(\eta,\xi)$.

For a pair $\zeta$ we define its {\it rotation number} $\rho(\zeta)\in[0,1]$ to be equal to the 
continued fraction $[r_0,r_1,\ldots]$ where $r_i=\chi({\cal R}^i\zeta)$. 
In this definition $1/\infty$ is understood as $0$, hence a rotation number is rational
if and only if only finitely many renormalizations of $\zeta$ are defined;
if $\chi(\zeta)=\infty$, $\rho(\zeta)=0$.
Thus defined, the rotation number of a commuting pair can be viewed as a rotation number in
the usual sense:
\begin{prop}
\label{rotation number}
The rotation number of the mapping $F_\zeta$ is equal to $\rho(\zeta)$.
\end{prop}

\noindent
There is an  advantage in defining $\rho(\zeta)$ using a sequence of heights in
removing the ambiguity in prescribing a continued fraction expansion to rational rotation numbers
in a renormalization-natural way.

\noindent
For $\rho=[r_0,r_1,\ldots]\in [0,1]$ let us set 
$$G(\rho)=[r_1,r_2,\ldots]=\left\{\frac{1}{\rho}\right\},$$
where $\{x\}$ denotes the fractional part of a real number $x$
($G$ is usually referred to as the {\it Gauss map}).
As follows from the definition,
$$\rho({\Ccal R}\zeta)=G(\rho(\zeta))$$ for a real commuting pair $\zeta$ with
$\rho(\zeta)\ne 0$.

The renormalization of the real commuting pair (\ref{real1}), associated
to some critical circle map $f$, is the rescaled pair
$(\wtl{{f}^{q_{m+2}}}|\wtl{{I}_{m+1}},\wtl{{f}^{q_{m+1}}}|\wtl{{I}_{m+2}})$.
Thus for a given critical circle map $f$ the renormalization operator
 recovers the (rescaled) sequence of the first return maps:
$$\{ (\wtl{f^{q_{i+1}}}|\wtl{I_i},\wtl{f^{q_{i}}}|\wtl{I_{i+1}})\}_{i=1}^{\infty}.$$

A {\it critical commuting pair} is a commuting pair $(\eta,\xi)$ whose maps  are real-analytic.
We shall also impose  a technical assumption that $\xi$  analytically extends to an interval $(a,b)\ni 0$ with 
$\xi(a,b)\supset [\eta(0),\xi(0)]$,
and has a single critical point $0$ in this interval. 
The space of critical commuting pairs modulo affine conjugacy  will be denoted by $\crit$; its subset consisting
of pairs $\zeta$ with $\chi(\zeta)=\infty$ will be denoted by ${\aaa S}_\infty$.
Renormalization is an injective transformation ${\cal R}:\crit\setminus\crit_\infty\to\crit$
(see \cite{Ya2}).

\section{Holomorphic extensions of critical commuting pairs}
\label{section: holomorphic extensions}
\noindent
In 1986 Eckmann and Epstein \cite{EE} introduced a space of critical commuting pairs now known as
the {\it Epstein class}. They showed that this class is invariant under the action of $\cR$, and
constructed the golden-mean fixed point of $\cR$ in this class using the methods of geometric complex
analysis. It was further shown by various people, such as Sullivan (in the unimodal case), \'Swia\c\negthinspace tek,
Herman, and Yoccoz, that the renormalizations of {\em any} $C^3$-smooth commuting pair with an irrational
rotation number converge to the Epstein class, at a geometric rate in the $C^2$-metric. 
Below, after some preliminaries, we define the Epstein class, and formulate these results more precisely.

\subsection{Carath\'eodory topology on a space of branched coverings.}
\label{triplets}
Consider the  collection $\bran$ of all triplets $(U,u,f)$, where $U\subset \Cbb C$ is a topological
disk different from the whole plane, $u\in U$, and $f:U\to \Cbb C$ is a three-fold analytic 
branched covering map, with the only branch point at $u$.
We will   topologize $\bran$ as follows (cf. \cite{McM1}).

Let $\{(U_n,u_n)\}$ be a sequence of open connected regions 
$U_n\subset \Cbb C$ with {\it marked points} $u_n\in U_n$.
Recall that this sequence {\it Carath\'eodory
converges} to a marked region $(U,u)$ if:
\begin{itemize}
\item $u_n\to u\in  U$, and
\item for any Hausdorff limit point $K$ of the sequence 
$\hat{\Cbb C}\setminus U_n$, $U$ is a component of 
$\hat{\Cbb C}\setminus K$.
\end{itemize}
For a simply connected $U\subset \Cbb C$ and $u\in U$ let $R_{(U,u)}:\Cbb D\to U$ denote
the inverse Riemann mapping with normalization $R_{(U,u)}(0)=u$, $R'_{(U,u)}(0)>0$.
By a classical result of Carath\'eodory, the Carath\'edory convergence
of simply-connected regions
$(U_n,u_n)\to (U,u)$ is equivalent to the locally uniform convergence of the inverse
Riemann mappings $R_{(U_n,u_n)}$  to $R_{(U,u)}$.

For positive numbers $\eps_1$, $\eps_2$, $\eps_3$ and compact subsets $K_1$ and $K_2$
of the open unit disk $\Cbb D$, let the neighborhood 
$\Ccal U_{\eps_1,\eps_2,\eps_3,K_1,K_2}(U,u,f)$ of an element $(U,u,f)\in \bran$
be the set of all $(V,v,g)\in\bran$, for which:
\begin{itemize}
\item $|u-v|<\eps_1$,
\item $\displaystyle \sup_{z\in K_1}|R_{(V,v)}(z)-R_{(U,u)}(z)|<\eps_2$,
\item and $R_{(U,u)}(K_2)\subset V$, and $\displaystyle \sup_{z\in R_{(U,u)}(K_2)}|f(z)-g(z)|<\eps_3$.
\end{itemize}
One verifies that the sets $\Ccal U_{\eps_1,\eps_2,\eps_3,K_1,K_2}(U,u,f)$ form
a base of a topology on $\bran$, which we will call {\it Carath\'eodory topology}. 
This topology is clearly Hausdorff, and  the convergence of a sequence $(U_n,u_n,f_n)$ to $(U,u,f)$
is equivalent to the Carath\'eodory convergence of the marked regions $(U_n,u_n)\to (U,u)$ as well
as a locally uniform convergence $f_n\to f$.

\noindent
\subsection{The Epstein class.}
An orientation preserving interval homeomorphism $g:I=[0,a]\to g(I)=J$ 
belongs to the {\it Epstein class $\Ccal E$} if it extends to an
analytic three-fold branched covering map of a topological disk
 $G\supset I$ 
onto the double-slit plane ${\Bbb C}_{\tl J}$, where $\tl J\supset \cl J$.
Any map $g$ in the Epstein class can be decomposed as
\begin{equation}
\label{epstein-decomposition}
g=Q_c\circ h,
\end{equation}
where $Q_c(z)=z^3+c$, and $h:I\to [0,b]$ is a univalent map 
$h:G\to \Delta(h)$ onto the complex plane with six slits,
which triple covers ${\Bbb C}_{\tl J}$ under the cubic map $Q_c(z)$.

For any $s\in (0,1)$, let us introduce a smaller class
 ${\Ccal E}_s\subset \Ccal E$ of Epstein mappings
$g:I=[0,a]\to J\subset \tl J$
for which both $|I|$ and $\dist(I,J)$ are $s^{-1}$-commensurable with $|J|$,
the length of each component of $\tl J\setminus J$
is at least $s|J|$, and $g'(a)>s$.
We will often refer to the space $\Ccal E$ as {\it the} Epstein class,
and to each ${\Ccal E}_s$ as {\it an} Epstein class.

We say that a commuting pair $(\eta,\xi)\in\crit$
belongs to the (an) Epstein class if both of its maps do.
Similarly, a critical circle map $f$ is Epstein if $\cR f$ is in the Epstein class.
It immediately follows from the definitions that:
\begin{lem} 
If a renormalizable commuting pair $\zeta$ is in the Epstein class, then the same is true for
 $\cR \zeta$.
\end{lem}

Let us make a note of an important compactness property of $\Ccal E_s$
\begin{lem}[Lemma 2.10 \cite{Ya2}]
\label{compactness}
Let $s\in (0,1)$. The collection of normalized maps $g\in {\Ccal E}_s$
with $I=[0,1]$,  with marked domains $(U,0)$ 
is sequentially compact with respect to Carath\'eodory
topology.
\end{lem}

\noindent
The importance of the Epstein class lies in the fact that all $C^1$-limit points of the sequence
$\{\cR^m(f)\}_{m=M}^\infty$ are in $\cE_s$ for a universal value of $s>0$. A more precise formulation
of this was proved in the recent work of de~Faria and de~Melo \cite{dFdM1}:
\begin{lem}
\label{real-bounds} 
There exists a universal constant $s>0$ such that the following holds.
Let $f\in C^r$, $(r\geq 3)$ be a critical circle map with an
 irrational rotation number.
Then the sequence of real commuting pairs 
$\cR^m(f)=(\wtl{{f}^{q_{m+1}}}|\wtl{{I}_{m}},\wtl{{f}^{q_{m}}}|\wtl{{I}_{m+1}})$
 is bounded in $C^{r-1}$-metric, and $C^{r-1}$-converges to $\cE_s$ at a geometric rate.
If $f\in C^\omega$, then the convergence is uniform on compact subsets of $\CC_{[0,1]}$.
\end{lem}

\noindent
In particular, for a critical circle map $f\in\cal E$ 
there exists $\sigma>0$ such that all its renormalizations are contained 
in $\cal E_\sigma$. Moreover, the constant $\sigma$ can be chosen
independent on $f$, after  skipping the first few renormalizations.

Finally, let us formulate an important statement about critical commuting pairs to be used further in the paper (Lemma 2.13, \cite{Ya2}): 
\begin{lem}[{\bf Parabolic Limits}]
\label{unique fixed point}
Let $\zeta=(\eta,\xi)\in \EE$ be a critical commuting pair with 
$\rho(\zeta)=0$, which appears as
 a limit of a sequence of critical commuting pairs $\{\zeta_n\}$
with $\rho(\zeta_n)\in \Cbb R\setminus \Cbb Q$.
 Then the map $\eta$ has a unique fixed point in
the interval $I_\eta$, which is necessarily parabolic, with
multiplier one.
\end{lem}

\noindent
A commuting pair $\zeta=(\eta,\xi)\in \EE$ will be called {\it parabolic} 
if the map $\eta$ has a unique fixed point in $I_\eta$, which has
a unit multiplier; this point will usually be denoted $p_\eta$.
Note, that
by virtue of its uniqueness, $p_\eta$ has to be globally attracting on one side
for the interval homeomorphism $\eta|_{I_\eta}$,
it is globally attracting on the other side under $\eta^{-1}$.

\subsection{Holomorphic commuting pairs}
\label{holomorphic pairs}
De~Faria \cite{dF1,dF2} introduced holomorphic commuting pairs to apply the
Sullivan's Riemann surface laminations argument to the renormalization of 
critical circle maps. They are suitably defined holomorphic extensions 
of critical commuting pairs which replace Douady-Hubbard polynomial-like maps \cite{DH2}.
A critical commuting pair $\zeta=(\eta|_{I_\eta},\xi|_{I_\xi})$  extends to a {\it
holomorphic commuting
pair} $\cal H$ if there exist four simply-connected $\RR$-symmetric domains $\Delta$, $D$, $U$, $V$ such that
\begin{itemize}
\item  $\bar D,\; \bar U,\; \bar V\subset \Delta$,
 $\bar U\cap \bar V=\{ 0\}$; the sets
 $U\setminus D$,  $V\setminus D$, $D\setminus U$, and $D\setminus V$ 
 are nonempty, 
connected, and simply-connected, $I_U=U\cap \RR\supset I_\eta$, $I_V=V\cap \RR\supset I_{\xi}$;
\item mappings $\eta:U\to (\Delta\setminus \RR)\cup\eta(I_U)$ and
 $\xi:V\to(\Delta\setminus \RR)\cup\xi(I_V)$ are onto and
univalent;
\item $\nu\equiv \eta\circ\xi:D\to (\Delta\setminus \RR)\cup{\nu(I_D)}$ is a three-fold 
branched covering with a unique critical point at zero, where 
$I_D=D\cap \RR$.
\end{itemize}

\begin{figure}
\begin{center}
\setlength{\unitlength}{0.00083333in}
{\renewcommand{\dashlinestretch}{30}
\begin{picture}(5942,3541)(0,-10)
\put(4184.020,3273.635){\arc{709.544}{1.0967}{2.1766}}
\put(3744.504,3545.229){\arc{1222.507}{1.1718}{1.9206}}
\put(3250.806,3257.465){\arc{807.040}{0.7894}{1.8054}}
\put(4024.350,2150.849){\arc{223.300}{4.2546}{6.2818}}
\put(4168.614,1878.160){\arc{548.467}{4.6115}{5.5534}}
\put(4578.804,1935.970){\arc{481.613}{3.6875}{4.8339}}
\put(4722.840,2258.898){\arc{284.444}{2.5106}{4.6713}}
\put(4760.529,2535.043){\arc{281.867}{1.8848}{4.3151}}
\put(4694.529,3017.400){\arc{699.884}{1.5151}{2.9826}}
\put(1043.483,2290.026){\arc{295.385}{5.5579}{7.5761}}
\put(1109.000,2648.000){\arc{509.902}{0.1974}{1.3734}}
\put(1485.961,2784.024){\arc{444.879}{0.8426}{2.1511}}
\put(1986.827,3003.491){\arc{1023.674}{1.3598}{2.3180}}
\put(973.812,1995.713){\arc{370.167}{5.3169}{6.8079}}
\put(1251.909,1807.853){\arc{295.308}{3.8418}{6.0078}}
\put(1487.306,1796.790){\arc{204.162}{3.6671}{5.6402}}
\put(1645.475,1853.164){\arc{153.257}{3.2047}{5.4857}}
\put(1750.563,1892.320){\arc{105.317}{3.3458}{7.2835}}
\put(1791.321,1802.107){\arc{95.036}{4.4501}{7.4728}}
\put(4184.020,232.365){\arc{709.543}{4.1066}{5.1865}}
\put(3744.504,-39.229){\arc{1222.508}{4.3626}{5.1114}}
\put(3250.806,248.535){\arc{807.040}{4.4778}{5.4938}}
\put(4024.350,1355.151){\arc{223.299}{0.0014}{2.0286}}
\put(4168.614,1627.840){\arc{548.468}{0.7298}{1.6717}}
\put(4578.804,1570.030){\arc{481.613}{1.4493}{2.5957}}
\put(4722.840,1247.102){\arc{284.444}{1.6119}{3.7725}}
\put(4760.529,970.957){\arc{281.867}{1.9681}{4.3984}}
\put(4694.529,488.600){\arc{699.884}{3.3006}{4.7681}}
\put(1043.483,1215.974){\arc{295.384}{4.9903}{7.0085}}
\put(1109.000,858.000){\arc{509.902}{4.9098}{6.0858}}
\put(1485.961,721.976){\arc{444.880}{4.1321}{5.4406}}
\put(1986.827,502.509){\arc{1023.674}{3.9652}{4.9233}}
\put(973.812,1510.287){\arc{370.167}{5.7585}{7.2495}}
\put(1251.909,1698.147){\arc{295.306}{0.2754}{2.4414}}
\put(1487.306,1709.210){\arc{204.161}{0.6430}{2.6161}}
\put(1645.475,1652.836){\arc{153.256}{0.7975}{3.0784}}
\put(1750.563,1613.680){\arc{105.315}{5.2829}{9.2205}}
\put(1791.321,1703.893){\arc{95.036}{5.0936}{8.1163}}
\put(1757.575,2687.476){\arc{704.349}{2.1767}{4.5683}}
\blacken\path(1602.388,2969.997)(1707.000,3036.000)(1583.635,3026.991)(1627.208,3009.746)(1602.388,2969.997)
\put(3560.882,2774.074){\arc{871.328}{4.6230}{6.7890}}
\blacken\path(3643.374,3231.838)(3522.000,3208.000)(3640.313,3171.916)(3605.891,3203.714)(3643.374,3231.838)
\put(2560.599,2479.135){\arc{494.410}{3.0275}{7.4090}}
\blacken\path(2301.005,2573.899)(2315.000,2451.000)(2360.486,2566.026)(2326.022,2534.274)(2301.005,2573.899)
\put(2524,1753){\ellipse{1190}{1190}}
\put(2971,1763){\ellipse{5926}{3510}}
\path(3155,2864)(2516,1757)
\path(2524,1753)(2089,2503)
\path(3155,642)(2516,1749)
\path(2524,1753)(2089,1003)
\thicklines
\path(1819,1753)(3654,1753)
\thinlines
\path(3972,2254)(3969,2256)(3965,2261)
	(3957,2268)(3947,2275)(3936,2283)
	(3924,2291)(3910,2297)(3895,2301)
	(3877,2303)(3856,2300)(3839,2294)
	(3825,2285)(3813,2275)(3803,2266)
	(3797,2258)(3794,2255)(3793,2254)
\path(3795,2254)(3794,2254)(3791,2255)
	(3782,2257)(3769,2260)(3754,2264)
	(3737,2266)(3721,2267)(3705,2266)
	(3690,2261)(3676,2253)(3668,2243)
	(3663,2232)(3661,2220)(3662,2208)
	(3663,2196)(3666,2184)(3669,2173)
	(3673,2163)(3675,2155)(3677,2150)(3678,2148)
\path(3683,2145)(3682,2145)(3677,2143)
	(3667,2140)(3653,2135)(3639,2128)
	(3625,2119)(3614,2107)(3606,2090)
	(3604,2071)(3607,2054)(3614,2039)
	(3622,2024)(3631,2012)(3638,2003)
	(3641,1999)(3642,1998)
\path(3648,1995)(3646,1994)(3643,1992)
	(3638,1989)(3632,1985)(3624,1979)
	(3617,1972)(3610,1964)(3604,1954)
	(3598,1941)(3595,1926)(3593,1907)
	(3594,1887)(3597,1867)(3602,1849)
	(3608,1831)(3615,1815)(3622,1800)
	(3629,1786)(3635,1773)(3640,1764)
	(3643,1759)(3645,1756)
\path(3972,1252)(3969,1250)(3965,1245)
	(3957,1238)(3947,1231)(3936,1223)
	(3924,1215)(3910,1209)(3895,1205)
	(3877,1203)(3856,1206)(3839,1212)
	(3825,1221)(3813,1231)(3803,1240)
	(3797,1248)(3794,1251)(3793,1252)
\path(3795,1252)(3794,1252)(3791,1251)
	(3782,1249)(3769,1246)(3754,1242)
	(3737,1240)(3721,1239)(3705,1240)
	(3690,1245)(3676,1253)(3668,1263)
	(3663,1274)(3661,1286)(3662,1298)
	(3663,1310)(3666,1322)(3669,1333)
	(3673,1343)(3675,1351)(3677,1356)(3678,1358)
\path(3683,1361)(3682,1361)(3677,1363)
	(3667,1366)(3653,1371)(3639,1378)
	(3625,1387)(3614,1399)(3606,1416)
	(3604,1435)(3607,1452)(3614,1467)
	(3622,1482)(3631,1494)(3638,1503)
	(3641,1507)(3642,1508)
\path(3648,1511)(3646,1512)(3643,1514)
	(3638,1517)(3632,1521)(3624,1527)
	(3617,1534)(3610,1542)(3604,1552)
	(3598,1565)(3595,1580)(3593,1599)
	(3594,1619)(3597,1639)(3602,1657)
	(3608,1675)(3615,1691)(3622,1706)
	(3629,1720)(3635,1733)(3640,1742)
	(3643,1747)(3645,1750)
\put(2499,2083){\makebox(0,0)[b]{$D$}}
\put(1469,2158){\makebox(0,0)[b]{$V$}}
\put(3524,2433){\makebox(0,0)[b]{$U$}}
\put(5244,2146){\makebox(0,0)[b]{$\Delta$}}
\put(1332,2728){\makebox(0,0)[rb]{$\xi$}}
\put(3972,3103){\makebox(0,0)[lb]{$\eta$}}
\put(2570,2863){\makebox(0,0)[b]{$\xi\circ\eta$}}
\put(2510,1483){\makebox(0,0)[b]{$0$}}
\path(484,1753)(5326,1752)
\end{picture}
}
\end{center}
\caption{A holomorphic commuting pair}
\label{holpair}
\end{figure}
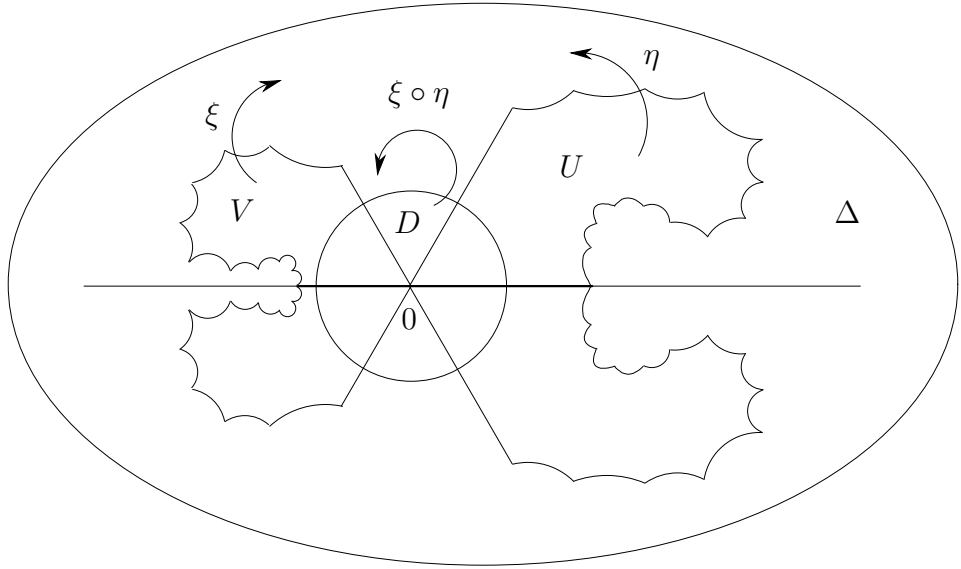

\noindent
We shall call $\zeta$ the {\it commuting pair underlying $\cH$}, and write $\zeta\equiv \zeta_\cH$.
The interval $I_\cH=[\eta(0),\xi(0)]$ will be called the {\it dynamical interval} of $\cH$.
The domain $D\cup U\cup V$ of a  holomorphic commuting pair $\cH$
will be  denoted $\Omega$ or $\Omega_{\cal H}$, the range will be denoted $\Delta$ or $\Delta_\cH$.
The closure of the  set of points whose 
orbit under $\cH$ is contained in $\Omega$
will be referred to as the {\it filled Julia set of $\cH$}, denoted $K({\cal H})$.
The {\it Julia set of $\cH$} is defined as  $J(\cH)=\partial K(\cH)$.

It is easy to see directly from the definition (cf. \cite{dF2}) that:
\begin{prop}Let $\zeta$ be a commuting pair with $\chi(\zeta)<\infty$.
Suppose $\zeta$ is a restriction of a holomorphic commuting pair $\cH$, that is
$\zeta=\zeta_{\cH}$. Then there exists a holomorphic commuting pair $\cG$ with range $\Delta_\cH$,
such that $\zeta_\cG=\cR\zeta$.
\end{prop}

\noindent
The {\it shadow} of the holomorphic commuting pair $\cH$ is the following piecewise-defined holomorphic
dynamical system:
$$S_{\cal H}(z) = \left\{ \begin{array}{l}
                   \eta(z),\; z\in U\\
                   \xi(z),\; z\in V\\
                   \xi\circ\eta(z),\; z\in D\setminus (U\cup V)
                 \end{array} \right. $$
As the next proposition shows one may think of the shadow of a holomorphic commuting
pair as an analogue of a cubic-like map:
\begin{prop}[Prop. II.4. \cite{dF2}]
\label{shadow-cover}
Given a holomorphic commuting pair $\cal H$ as above, consider its shadow 
$S_{\cal H}$.
Let $I=\Omega \cap \Bbb R$, and $X=I\cup S_{\cal H}^{-1}(I)$. Then:
\begin{itemize}
\item The restriction of $S_{\cal H}$ to $\Omega\setminus X$ 
is a regular three fold
covering onto $\Delta\setminus \Bbb R$.
\item $S_{\cal H}$ and $\cal H$ share the same orbits as sets.
\end{itemize}
\end{prop}

\noindent
We will say that two holomorphic commuting pairs
${\cal H}:\Omega_{\cal H}\to \Delta_{\cal H}$ and 
${\cal G}:\Omega_{\cal G}\to \Delta_{\cal G}$ are conjugate if there 
is a homeomorphism $h:\Delta_{\cal G}\to \Delta_{\cal H}$ such that
$$S_{\cal G}=h^{-1}\circ S_{\cal H}\circ h.$$
In this case we will simply write ${\cal G}=h^{-1}\circ {\cal H}\circ h.$

\medskip
\subsection{Complex {\it a priori} bounds}
\label{sec bounds}
We shall denote by $\hol$ the space of holomorphic commuting pairs
$\Ccal H:\Omega\to \Delta$ whose underlying real commuting pair $(\eta,\xi)$
is in the Epstein class. 
In this case both maps $\eta$ and $\xi$ extend to triple branched coverings 
$\hat{\eta}:\hat{U}\to \Delta\cap {\Bbb C}_{\eta(J_\eta)}$ and 
$\hat{\xi}:\hat{V}\to\Delta\cap {\Bbb C}_{\xi(J_\xi)}$ respectively.
We will turn $\hol$ into a topological space by identifying it with
a subset of 
$\bran\times\bran$ by $\Ccal H\mapsto (\hat{U},0,\hat{\eta})\times(\hat{V},0,\hat{\xi})$
(cf. \secref{triplets}).

\noindent
We say that a real commuting pair $(\eta,\xi)$ with
an irrational rotation number has
{\it complex {\rm a priori} bounds}, if all its renormalizations extend to 
holomorphic commuting pairs with {\it bounded moduli}:
$$\mod(\Delta\setminus\Omega)>\mu>0.$$
For $\mu\in(0,1)$ let $\hol(\mu)$ denote the space of holomorphic commuting pairs
${\cH}:\Omega_{{\cH}}\to \Delta_{\cH}$, with $\mod (\Delta_{{\cH}}\setminus\Omega_{{\cH}})>\mu$,
$\min(|I_\eta|,|I_\xi|)>\mu$ and $\diam(\Delta_{\Ccal H})<1/\mu$.

\begin{lem}[Lemma 2.15 \cite{Ya2}]
\label{bounds compactness}
For each $\mu\in(0,1)$ the space $\hol(\mu)$ is sequentially pre-compact, with every limit point
contained in $\hol(\mu/2)$.
\end{lem}

\noindent
The existense of complex {\it a priori} bounds is a key analytic issue 
 of renormalization theory. In the case of critical circle maps it is settled
by the following theorem:

\begin{thm}
\label{complex bounds}There exist universal constants $\mu>0$ and $K>1$ such that
the following holds. 
Let $\zeta\in\crit$ be a critical commuting pair with an irrational rotation number. 
Then there exists $N=N(\zeta)$ such that for all
$n\geq N$ the  commuting pair $\cR^n\zeta$ extends to a holomorphic commuting pair
$\cH_n:\Omega_n\to\Delta_n$ in ${\aaa H}(\mu)$.
The range $\Delta_n$ is a Euclidean disk of radius at most $K$, and the regions $\Omega_n\cap(\pm \HH)$
are $K$-quasidisks.
\end{thm}

\begin{rem}
\label{rem bounds}
We first proved this theorem in \cite{Ya1} for commuting pairs 
$\zeta$ in an Epstein class $\cE_s$, in which case
 $N=N(s)$. Our proof was later adapted by de~Faria and de~Melo \cite{dFdM2} 
to the case of a non-Epstein critical commuting pair. In the general case, in a Carath{\'e}odory
compact family of critical commuting pairs, the number $N$ can be chosen uniformly.
\end{rem}

\noindent
Let $\zeta$ be at least $n$ times renormalizable
 critical commuting pair. For the lack of a better term, let us say that
the pair of numbers $\tau_n(\zeta)=(r_{n-1},r_{n-2})$ forms the {\it history} of
the pair $\Ccal R^n\zeta$.
Based on the above theorem and a detailed analysis of the shapes of the domains $\Omega_n$ we 
proved the following in \cite{Ya1}:

\begin{thm}[\cite{Ya1}]
\label{qc-conjugacy}
There exists a universal constant $K_1>1$ such that the following holds.
Let $\zeta_1=(\eta_1,\xi_1)$ and $\zeta_2=(\eta_2,\xi_2)$ be two  critical 
commuting pairs with irrational rotation numbers. 
Let $n>max(N(\zeta_1),N(\zeta_2))+1$ as above.
Assume that the $n$-th renormalizations of $\zeta_1$, $\zeta_2$ have the same rotation number
and the same history. Then their holomorphic commuting pair extensions $\cH_n^1$, $\cH^2_n$ are $K_1-$quasiconformally conjugate. The conjugating map is conformal on the filled Julia
set.
\end{thm}
\noindent
For commuting pairs of  the type bounded by $B$ this theorem was first proved by 
de~Faria \cite{dF1,dF2},
with ``$K_1$" depending on the value of $B$.
The proof of the above theorem in the case of an unbounded type rotation number
requires an analysis of the shape of the domain of the holomorphic pair. 
In what follows, let $f$ be an analytic critical circle mapping with $\rho(f)\in\RR\setminus\QQ$,
 fix a sufficiently large $n$ and let $\cH:\Omega\to\Delta$ be the holomorphic 
pair extension of $\cR^n f$ guaranteed by \thmref{complex bounds}.

Consider the inverse orbit:
\begin{equation}
\label{JJ-orbit}
J_0\equiv f^{q_{n+1}}(I_n),J_{-1}\equiv f^{q_{n+1}-1}(I_n),\ldots,J_{-(q_{n+1}-1)}\equiv f(I_n),
\end{equation}
and  the corresponding inverse orbit for the domain $\Delta_0\equiv
\Delta\cap \Bbb H$:
\begin{equation}
\label{delta-orbit}
\Delta_0,\Delta_{-1},\ldots,\Delta_{q_{n+1}-(-1)}\equiv f(U)\cap \Bbb H.
\end{equation}
Consider the consecutive returns of the orbit (\ref{JJ-orbit})
to $I_{m-1}$ before the first return to $I_m$,
\begin{equation}
\label{returns}
J_{-q_m},J_{-2q_m},\ldots,J_{-l_m q_m}.
\end{equation}
 Consider the curve segment 
$\gamma^m\subset f^{-q_m}([f^{q_{m-1}-q_m}(0),
f^{q_{m-1}}(0)])\cap \Bbb H$, $\cl \gamma^m\ni f^{q_{m-1}-q_{m}}(0)$, and let

\begin{equation}
\label{gamma-orbit}
\gamma^m\equiv\gamma^m_{-1},\gamma^m_{-2},\ldots,\gamma^m_{-l_m}
\end{equation}
 be the corresponding
inverse orbit of the curve segment $\gamma^m$ under $f^{q_m}$. Let the
curve $\Gamma^m$ be the union of the segments $\gamma^m_i$, for $i=-1,\ldots,
-l_m$.

\begin{figure}
\begin{center}
\input figpic/pinch.eepic
\end{center}
\caption{}
\label{pinchfig}
\end{figure}

\begin{lem}
\label{cut-out}
There exists a topological  disk $\hat D\subset  D([f^{sq_n-q_{n+1}}(0),
f^{-q_{n-1}-sq_n}(0)])$ commensurable with $I_{n-1}$,
such that 
\begin{itemize}
\item The domain $U\setminus \hat D$ is a $K$-quasidisk for some fixed $K$.
\item  The intersection $\partial U\cap \hat D\cap \Bbb H$ is contained
in the curve $f^{-(q_{n-1}-1)}(\Gamma^n)$, where $\Gamma^n$ is  as above.
\end{itemize}
\end{lem}

\noindent
In conclusion, let us say that a holomorphic pair $\cH$ is {\it $K$-bounded} if the conclusions of \thmref{complex bounds}
and \lemref{cut-out} hold for $\cH$ with this value of $K$ and $\mu=1/K$, and each of the maps $\xi,$ $\eta$, and $\nu$ 
constituting $\cH$ is a composition of the cubic map $z\mapsto z^3$ and a conformal diffeomorphism with $K$-bounded distortion.

\section{McMullen's rigidity result}
We are going to briefly review here some of the results of McMullen on geometric limits
in dynamics, as discussed in Chapter 9 of \cite{McM2}. The definitions are at times
quite technical, therefore we will not attempt to go into every detail, and will
instead try to outline the main points relevant to our analysis. McMullen gives a general
definition of a {\it holomorphic dynamical system} $\cF$ in $\hat\CC$ as a countable collection of analytic
hypersurfaces in $\hat \CC\times\hat\CC$ which define the holomorphic relationships forming $\cF$.
He endows the space of holomorphic dynamical systems with the Hausdorff topology,
the convergence in this topology is referred to as the {\it geometric convergence}.
Following \cite{dFdM2},
we will consider the holomorphic dynamical system $\cF=\cF(\cH)$ generated by the complete dynamics
of a holomorphic pair $\cH$. Each hypersurface in $\cF(\cH)$ is
the graph of a relation $(S_{\cH})^i(z)=(S_\cH)^j(z)$ for some $i,j\in\ZZ$.

Recall that a line field $\mu$ on $\hat\CC$ is {\it parabolic} if $\mu=A^*(dz/d\bar z)$ where 
$A\in\operatorname{Aut}(\hat{\CC})$. A dynamical system $\cF$ is {\it nonlinear} if no
parabolic line field is left invariant by $\cF$. 
It is {\it twisting} if any holomorphic dynamical system $\cF_1$ which
is quasiconformally conjugate to $\cF$ is non-linear. 
A quantitative measure of non-linearity may be introduced as follows.
Let $$\sg(z)\abs{dz}=\frac{2\abs{dz}}{1+\abs{z}^2}$$
be the spherical metric and denote $B(x,r)$ the spherical ball of radius
$r$ centered at $x$. Let $B_1,\ldots,B_k$ be any collection of disjoint spherical balls,
and let $f_i:B_i\to\hat\CC$ be a univalent branch in $\cF$ which extends to a ball with twice the
radius. 
The nonlinearity $\nu(\cF)$ of $\cF$ is defined as
\begin{equation}
\label{nonlinearity}
\nu(\cF)=\inf_{\text{parabolic }\mu}\sup_{\text{collections} f_1,\ldots,f_k}
\sum_{i=1}^k\int_{B_i}|{\mu-(f_i)^*\mu}|\sigma^2(z)dxdy.
\end{equation}
Our definition varies slightly form the definition of \cite{McM2} who considers supremum over
the collections of maps with only a single element. 
It is not difficult to see that this change does not affect the conclusions. As an example, the reader is invited
to observe how the proof of Lemma 9.12 of \cite{McM2} will look with our definition (this lemma is a
key element of the proof of the Inflexibility Theorem below). Instead of passing to the limit of a subsequence
$\cF^{\text{sat}}_n\to\cG$ as McMullen does, we consider a blow up $v_n$ of $v$, which is nearly invariant
under a uniformly twisting system $\cG_n$ in the measurable sense. This leads to a contradiction in the
same way as before.

For a holomorphic dynamical system $\cF$ and $K>1$ set
$$
\nu^K(\cF)=\inf_\phi\nu(\phi_*(\cF)),
$$
where infinum is taken over $K$-quasiconformal mappings
$\phi:\hat\CC\to\hat\CC$, $\phi(0)=0$, $\phi(1)=1$,
$\phi(\infty)=\infty$. 
A collection of holomorphic dynamical systems
$\cF_\alpha$ is {\it uniformly twisting} 
if for any $K>1$ 
$$\inf_\alpha \nu^K(\cF_\alpha)>0.$$
A holomorphic commuting pair $\cH$ is uniformly twisting if any collection of blow-ups 
of $\cF=\cF(\cH)$ at the points of $J(\cH)$ is uniformly twisting. More precisely,
let $\Lambda(\cH)$ be the set of all pairs
$(x,r)\in\CC\times\RR_+$ such that $x\in J(\cH)$, $0\le r\le 1$. Let
$$
\nu^K(\cF,\Lambda(\cH))=\inf_{\omega\in\operatorname{ch}(\Lambda(\cH))}
\nu^K(T^\omega_*(\cF)),
$$
where $\operatorname{ch}(\Lambda(\cH))$ is the convex hull of $\Lambda$ in
$\CC\times\RR_+$ and $T^\omega$, $\omega=(x,r)$ is the fractional linear
transformation moving the ball $B(x,r)$ to $B(0,1)$.
Then $\cH$ is called {\it uniformly twisting}
if $\nu^K(\cF,\Lambda(\cH))>0$ for all $K>1$. 

\noindent
The relevance of the above definitions lies in the rigidity result below; to formulate it
we need another definition.
A point $x\in J(\cH)$ is a {\it $\beta$-deep point} if there exists $\beta>0$ such that for
every spherical ball $B(x,r)$ with $r$ sufficiently small, 
the largest ball contained in $B(x,r)\setminus J(\cH)$ has radius $s(r)\leq r^{1+\beta}$.

\begin{thm}[{\bf Dynamical inflexibility} \cite{McM2}]
\label{thm-inflex}
  Let $(\cF(\cH),\Lambda(\cH))$ be uniformly twisting and let
  $\phi:\hat\CC\to\hat\CC$ be a $K$-quasiconformal conjugacy between holomorphic pairs $\cH$
  and $\cG$. Then for any $\beta$-deep point $x$ of $J(\cH)$, $\phi$
  is $C^{1+\alpha}$ conformal at $x$ and constant $\alpha$ depends only on
  $K$, $\dl$ and $\nu^K(\cF(\cH),\Lambda(\cH))$.
\end{thm}

\noindent
de~Faria and de~Melo showed (see \cite{dFdM2}):
\begin{thm}
\label{twisting bounded type}
Let $\cH$ be a holomorphic commuting pair with an irrational rotation number $\rho$.
There exists $\beta>0$ such that the following holds:
\begin{itemize}
\item[(I)] the critical point $0$ is a $\beta$-deep point of $J(\cH)$;
\item[(II)] if $\rho$ is of bounded type, then $(\cF(\cH),J(\cH))$ is uniformly twisting.
\end{itemize}
\end{thm}

\noindent
To prove  part (II) of \thmref{twisting bounded type} de~Faria and de~Melo show:

\begin{lem}[\cite{dFdM2}]
\label{dFdM-lemma}
For every $A\in\NN$
there
exist constants $K>1$, $\delta>0$, $\mu>0$, $C>0$ such that
 for any holomorphic pair $\cH$ with an irrational rotation number of a type
bounded by $A$, 
 for every point $z\in J(\cH)$,
and every $0<r<\delta$, the disk $B(z,r)$ contains a holomorphic pair $\cH_1:\Omega_1\to\Delta_1\in\cF(\cH)$
with the complex {\it a priori} bound $\mu$,
such that $\diam(\Delta_1)>Cr$, and $\Delta_1$ is a $K$-quasidisk. 
\end{lem}

\noindent
They then use the following fact:
\begin{prop}[\cite{dFdM2}]
\label{dFdM-prop}
In the notation of the previous lemma, 
there exists $\upsilon>0$ depending only on the values of $K$, $\mu$, and $C$,
such that the nonlinearity of $S_{\cH_1}$ is bounded below by $\upsilon\cdot\area{\Delta_1}$
\end{prop}

\noindent
Taking the branches of $S_{\cH_1}$ to be the collection of maps in (\ref{nonlinearity}), we see that
the nonlinearity of $\cF$ in the disk $B(z,r)$ is bounded from below by $\operatorname{const}\cdot\operatorname{area}(B(z,r))$. After rescaling $B(z,r)$ to $B(0,1)$, we see that the nonlinearity is bounded below 
by a uniform constant, and hence $\cH$ is uniformly twisting.

\noindent
In this paper we will demonstrate:
\begin{thm}
\label{uniform twisting}
There exists $\nu_0>0$ such that for every holomorphic pair $\cH$ with an irrational
$\rho$, $(\cF(\cH),J(\cH))$ is uniformly twisting, with
$$\nu^{K_1}(\cF(\cH),\Lambda(\cH))\geq \nu_0,$$
where $K_1$ is as in \thmref{qc-conjugacy}

\end{thm}

\noindent
To understand the difficulty involved in proving uniform twisting in the case of a 
rotation number of unbounded type, we need first to discuss the local theory
of parabolic perturbations. After a brief discussion, in the next section we will outline
the idea of the proof of \thmref{uniform twisting}.

\section{Parabolic  maps and their perturbations}
\label{parabolic}

\subsection{General facts}
      We begin with a brief review of the theory of parabolic bifurcations,
as applied in particular to an interval map in the Epstein class. For a more
comprehensive exposition the reader is referred to \cite{Do},
supporting technical details may be found in \cite{Sh}.
Fix a  map $\eta_0\in\Ccal E$ having a parabolic fixed point $p$
with unit multiplier.

\begin{thm}[Fatou Coordinates]
\label{Fatou-coord}
There exist topological discs $U^A$ and $U^R$, called {\it attracting}
and {\it repelling petals}, whose union is a 
punctured neighborhood of the parabolic periodic point $p$
such that
$$\eta_0(\bar U^A)\subset U^A\bigcup \{p\},\text{ and }
\bigcap_{k=0}^\infty \eta_0^{k}(\bar{ U}^A)=\{p\},$$

$$\eta_0(\bar U^R)\subset U^R\bigcup \{p\},\text{ and }
\bigcap_{k=0}^\infty \eta_0^{-k}(\bar{ U}^R)=\{p\},$$
where $\eta_0^{-1}$ is the univalent branch fixing $\zeta$.

\noindent
Moreover, there exist injective analytic maps
$$\Phi^A:U^A\to{\Bbb C}\text{ and }\Phi^R:U^R\to \Bbb C,$$ 
unique up to post-composition by translations, such that
$$\Phi^A(\eta_0(z))=\Phi^A(z)+1\text{ and }\Phi^R(\eta_0(z))=\Phi^R(z)-1.$$
The Riemann surfaces $C^A=U^A/\eta_0$ and $C^R=U^R/\eta_0$ are
 conformally equivalent to the cylinder $\Bbb C/\Bbb Z$.

The coordinate change $\Phi^A(z)=-1/(z-p)+o(1/(z-p))$ and similarly for $\Phi^R$.
\end{thm}

\noindent
We denote $\pi_A:U^A\to C^A$ and $\pi_R:U^R\to C^R$ the natural projections.
The quotients $C^A$ and $C^R$ are customarily referred to as {\it \'Ecalle-
Voronin cylinders}; 
 we will find it useful to regard these as Riemann
spheres with distinguished points $+,-$ filling in the punctures.
The real axis projects to the {\it natural equators} $E^A\subset C^A$
and $E^R\subset C^R$.
Any conformal {\it transit homeomorphism} $\tau:C^A\to C^R$ fixing the ends $+,-$
is a translation in suitable coordinates.  Lifiting it produces a map
 $\bar \tau:U^A\to \Cbb C$ satisfying 
$$\tau\circ\pi_A=\pi_R\circ\bar\tau.$$
We will sometimes write $\tau\equiv \tau_\theta$, and $\bar\tau=\bar\tau_\theta$, where
$$\Phi^R\circ\bar\tau\circ(\Phi^A)^{-1}(z)\equiv z+\theta \mod\ZZ.$$
The return map from $U^R$
to $U^A$ descends to a well-defined analytic transformation 
$${\cal E}:{\cal W}\to C^A$$
(the {\it {\'E}calle-Voronin map})
where ${\cal W}$ is an open subset of $C^R$. It is easy to
see that the ends of $C^R$ belong to different components of 
${\cal W}$. The choice of a conformal
{\em transit isomorphism} $$\Theta:C^A\to C^R$$ respecting these ends
determines an analytic dynamical system
$${\cal F}_\Theta=\Theta\circ{\cal E}:{\cal W}\to C^R$$
with fixed points at $\pm$.  The product of the corresponding
eigenvalues $\varrho^\pm_{\Theta}$ is clearly independent of $\Theta$, and by the Schwarz Lemma
is a number greater than one.

\noindent 
Suppose
for an analytic  map $\eta$ in a sufficiently small neighborhood of $\eta_0$
the parabolic point splits into a complex conjugate pair of repelling
fixed points $p_\eta\in \Bbb H$ and $\bar p_\eta$ with multipliers
$\lambda_\eta^\pm=e^{2\pi i\pm\alpha(\eta)}$. In this situation one may still speak of 
attracting and repelling petals:

\begin{lem}[{\bf Douady Coordinates}]
\label{douady-coord}
Let $V\subset \CC$ be a domain containing $I_{\eta_0}$.
There exists a Carath{\' e}odory neighborhood $U(\eta_0)$ of the map $\eta_0$
in  the domain $V$ such that the following holds.
For any $\eta\in U(\eta_0)$ with $|\arg \alpha(\eta)|<\pi/4$,
there exist topological discs $U_\eta^A$ and $U_\eta^R$ whose union is a 
neighborhood of $p$, and injective analytic maps
$$\Phi^A_\eta:U^A\to {\Bbb C}\text{ and }\Phi_\eta^R:U^R_f\to{\Bbb C}$$
unique up to post-composition by translations, such that
$$\Phi^A_\eta(\eta(z))=\Phi^A_\eta(z)+1
\text{ and }\Phi_\eta^R(\eta(z))=\Phi_\eta^R(z)+1.$$
The quotients $C^A_\eta=U^A_\eta/\eta$ and $C^R_\eta=U^R_\eta/\eta$
 are Riemann surfaces
conformally equivalent to ${\Bbb C}/{\Bbb Z}$. 
\end{lem}

\noindent
Let us note:
\begin{prop}
There exists an open neighborhood $\cW(\eta_0)$ of $\eta_0$ in the 
Carath\'eodory topology in the domain $V$, such that for every $\eta\in\cW(\eta_0)$ as above,
the condition on the eigenvalues of the repelling fixed points is automatically
satisfied.
\end{prop}

An arbitrary choice of real basepoints $a\in U^A$ and $r\in U^R$ enables us
to specify the Fatou and Douady coordinates uniquely, by requiring that
$\Phi^A(a)=\Phi^A_\eta(a)=0$, and $\Phi^R(r)=\Phi^R_\eta(r)=0$.
The following fundamental theorem first appeared in 
\cite{orsay-notes}:

\begin{thm}
\label{continuity-fatou}
With these normalizations the maps $\Phi_\eta^A$, $\Phi^R_\eta$ depend continuously on $\eta$
with respect to the compact-open topology, and
$$\Phi^A_\eta\to\Phi^A\text{ and }\Phi^R_\eta\to \Phi^R$$
uniformly on compact subsets of $U^A$ and $U^R$ respectively.

\noindent
Moreover, select the smallest $n(\eta)\in\Cbb N$ for which 
$\eta^{n(\eta)}(a)\geq r$. Then 
$$\eta^{n(\eta)}(z)=(\Phi^R_\eta)^{-1}\circ T_{\theta(\eta)+K}\circ \Phi^A_\eta$$
wherever both sides are defined. In this formula $T_a(z)$ denotes the
translation $z\mapsto z+a$, $\theta(\eta)\in [0,1)$ is given by 
$$\displaystyle\theta(\eta)=1/\alpha(\eta)+
\underset{\alpha(\eta)\to\infty}{o(1)}\mod 1,$$ and the real
constant $K$ is determined by the choice of the basepoints $a$, $r$.
Thus for a sequence $\{\eta_k\}\subset U(\eta)$ converging to $\eta$, 
the iterates $\eta_k^{n({\eta_k})}$ converge locally uniformly if and only
if there is a convergence $\theta(\eta)\to \theta$, and the limit 
in this case is a certain lift of the transit homeomorphism $\tau_\theta$
for the parabolic map $\eta_0$.
\end{thm}

\begin{figure}
\input{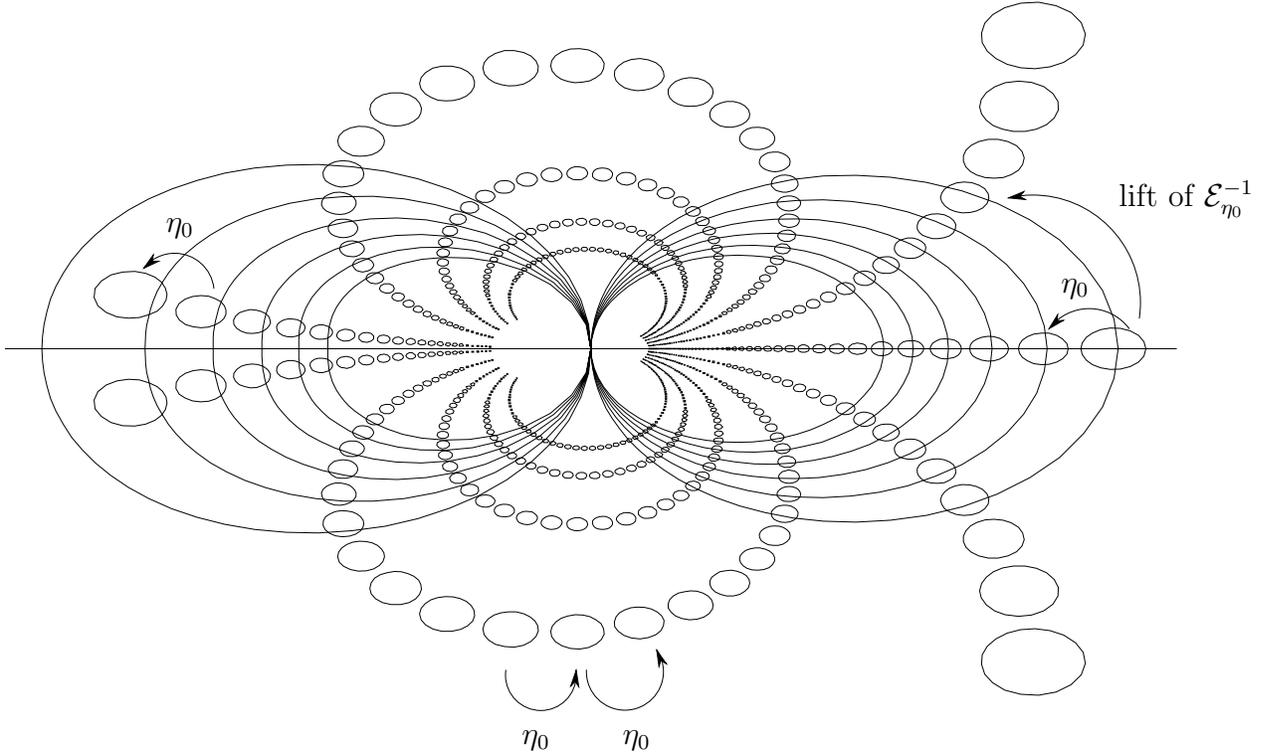}
\caption{Images of a disk under the parabolic dynamics}
\oplabel{f:parab}
\end{figure}

\subsection{An explanation of the method of the proof of \thmref{uniform twisting}}
Let us begin by explaining the difference between the bounded and unbounded cases. 
To understand what happens in the case of an unbounded type, consider a sequence of renormalizable 
holomorphic pairs
$\cH_n\to\cH$  with $\rho(\cH)=0$, such that $\cR(\cH_n)$ also converge to a holomorphic pair $\cH'$. 
In the Fatou coordinates of $\cH$, the preimages of of $\cH'$ form a grid $\Lambda$, which up to a 
bounded distortion 
is produced by the translations $z\mapsto z+1$, $z\mapsto z+\log(\cE'(+))/2\pi i$ (cf. Figure \ref{f:parab}).

\noindent
Since the Fatou coordinates have the order $1/(z-p)$, the size of the largest holomorphic commuting pair of the
grid in the disk $D_r(p)$ is of the order $r^2$. Thus in this case we cannot get the uniform twisting condition
by finding a single holomorphic pair commensurable with the disk, as in \lemref{dFdM-lemma}. This is precisely where the
argument of \cite{dFdM2} fails in the case of unbounded type.

The problem is resolved in this paper using the following simple idea. In any disk in $\CC$, the domains of the
holomorphic pairs of the grid $\Lambda$ take up a universally bounded proportion of the area. 
Therefore, the same is true
in each of the disks $D_r(p)$ for small enough $r$. We will thus use {\it all} of the copies of $\cH'$ 
rather than just one of them to prove that $\cH$ is uniformly twisting.

\def\Eeta{\cE_{\eta}}
\def\tDl{\tilde\Delta}
\def\lletter{l}

\section{Key technical lemma}
\label{s:proof}


\begin{figure}
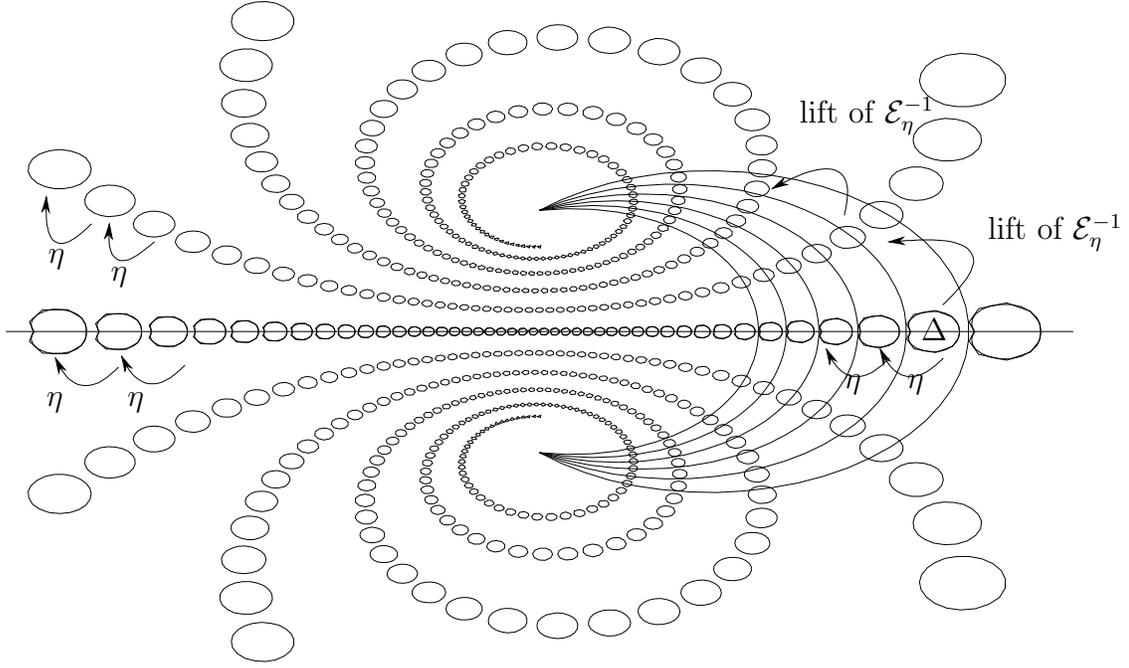

{
\def\phi{\eta}
\input pic/parpic.eepic

}
\caption{Images of disk $\Delta$ under the full dynamics of $\eta$}
\label{f:diskDelta}
\end{figure}

\noindent
For a critical commuting pair $(\eta,\xi)$ set
$I=I_\eta\setminus\eta(I_\eta)$.  Notice that the largest $r$ for
which $\eta^r(I)\subset{}I_\eta$ is equal to the height $\chi(\zeta)$.

\begin{lem}
  \label{l:grade_for_almost_parabolic_map}
  For every $K>1$ there exists $0<C<1$ such that the
  following holds.  Suppose $\zeta$ is a critical commuting pair with
  an irrational rotation number.  Then there exists $M$ such that for
  every $m>M$ denoting $(\eta,\xi)=\cR^m\zeta$ and 
  $r=\chi(\eta,\xi)$, we have the following.
  
  Let $\Delta$ be a disk $K$-commensurable with one of the intervals
  $\eta^{\varkappa}(I)$, $0\le \varkappa\le r$, such that $\Delta$ is symmetric with
  respect to the real axis, $\Delta\cap\eta^{\varkappa}(I)\ne\varnothing$ and
  $(\Delta\setminus\eta^{\varkappa}(I))\cap\RR=\varnothing$, $\mod(\Delta)\ge
  1/K>0$.  Then for any $z_0\in\eta^i(I)$, $0\le{i}\le{}r$ for any
  $\lletter\in[\abs{\eta^i(I)},1]$ there exists a pointed area $(U,y)$
  such that $\abs{z_0-y}\sim{\lletter}$, $\diam U\sim{\lletter}$ and
  $$
  \sum_{\tilde\Delta\subset U} \area{\tilde\Delta} \ge C\lletter^2,
  $$
  where $\phi\tilde\Delta=\Delta$ for some $\phi$ from the
  complete dynamics generated by $\eta$.
\end{lem}

\begin{proof}
  It is instructive to have a look at Figure~\ref{f:diskDelta}, to see
  how $\Delta$ is moved around under $\eta$ and $\eta^{-1}$. Notice
  that the mapping $\Eeta$ is some pertubation of Ecalle-Voronin
  mapping $\cE$. In fact, $\Eeta$ is the return map ${R}_\eta$ to a
  strip connecting $z_-$ to $z_+$. To describe the picture
  mathematically, we have to pass to the Douady coordinates.  First,
  Theorem~\ref{complex bounds} guarantees that there exists $M>0$ such
  that for all $m>M$ commuting pair $\cR^m\zeta$ extends to a
  holomorphic commuting pair $\cH_m$ with universal complex {\it a
    priori} bounds.  Hence we can apply the following fact, which
  follows from \lemref{douady-coord} and \thmref{continuity-fatou}.

  \begin{figure}
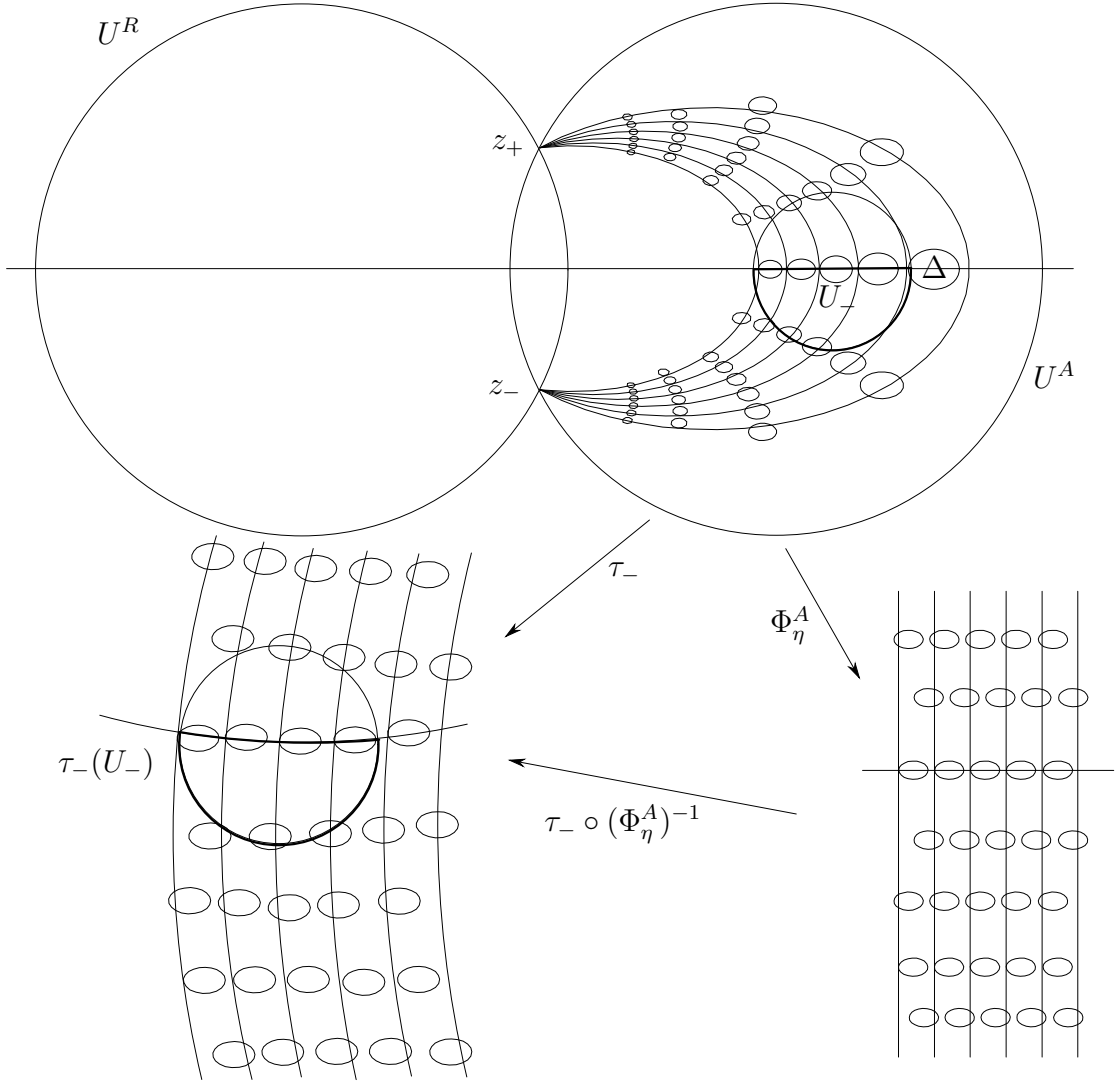

    \centering
    {
      \parindent=-1cm
      \input figpic/tauphi.eepic
      
    }
    \caption{Estimate for the area of discs $\tDl$, covered by $U_-$}
    \label{f:phi_tau-}
  \end{figure}

  There exist universal $\dl>0$, $N>0$, such that for all
  $r>N$ the following statements hold. Denote $\Phi^A_\eta$ and $\Phi^R_\eta$
the attracting and repelling Douady coordinates of $\eta$, defined in 
  $U^A$ and $U^R$ respectively, where $U^A$ and $U^R$ are real-symmetric disks
of radius $\delta$ whose boundaries contain both of the points $z_+$, $z_-$.


   For $\tau_-(z)=z_- -1/(z-z_-)$ and $\tau_+(z)=z_+ -1/(z-z_+)$ we
  have 
  $$\frac1K\le(\tau_-\circ(\Phi^A_\eta)^{-1})'(w)\le{}K\text{ for }w\in\Phi^A_\eta(U^A\cap(-\HH)),$$
  $$\frac1K\le(\tau_+\circ(\Phi^A_\eta)^{-1})'(w)\le{}K\text{ for }w\in\Phi^A_\eta(U^A\cap\HH),$$
  $$\frac1K\le(\tau_-\circ(\Phi^R_\eta)^{-1})'(w)\le{}K\text{ for }w\in\Phi^R_\eta(U^R\cap(-\HH)),$$
  $$\frac1K\le(\tau_+\circ(\Phi^R_\eta)^{-1})'(w)\le{}K\text{ for }w\in\Phi^R_\eta(U^R\cap\HH).$$
  
  Notice that if $r\le{}N$, then all images of $\Delta$ by $\eta$ and
  its inverse are commensurable by the Koebe Distortion Theorem, hence
  the statement is evident. Therefore we assume $r>N$. Since $\dl$ is
  universal, we can restrict ourselves to
  $z_0\in(U^A\cup{}U^R)\cap\eta^i(I)$, and $\lletter\le\dl$. By
  symmetry consider the case $z_0\in{}U^A$ only.  Given
  $\lletter\in[\abs{\eta^i(I)},\dl/2]$ one can find real
  $y\in\eta^j(\Delta\cap\RR)$, such that $\abs{y-z_0}\sim\lletter$ and
  the circle $U$ of radius $\lletter$ about $y$ lies strictly inside
  $U^A$.
  
  Let us study the intersections of $U_-=U\cap(-\HH)$ with $\tDl$.
  Firstly, $U_-$ intersects halfs of $\tDl$, which are the images and
  preimages of $\Delta$ by $\eta$. However, $\eta$ also has a
  branching point at 0, which gives rise to $\tDl$ lying in $-\HH$ or
  $\HH$ (see Figure~\ref{f:diskDelta}). In fact, it follows from
  \lemref{douady-coord}, that the set $\{\Phi^A_\eta(\tDl)\}$ is a
  bounded distortion image of the lattice
  $\sqcup_{n,m\in\ZZ}D_{1/4}(n+im)$ (see
  Figure~\ref{f:phi_tau-}). By the above estimates, the set
  $\{\tau_-\circ(\Phi^A_\eta)^{-1}\circ\Phi^A_\eta(\tDl)\}=\{\tau_-(\tDl)\}$
  forms a perturbed lattice in $\tau_-(U^A\cap(-\HH))$. Consider
  \begin{equation}
    \label{e:sumtDl1}
    \sum_{\tDl} \area{\tDl\cap U_-}=\sum_{\tDl}\iint_{\tau_-(\tDl\cap
      U_-)}J_{\tau_-^{-1}}(w)dudv,
  \end{equation}
  where $w=u+iv$, $dudv$ is an element of area, $J_{\tau_-^{-1}}(w)$
  is the Jacobian of $\tau_-^{-1}$ at $w$. Since $\tau_-^2$ is the identity,
  $$
  \tau_-^{-1}(w)=z_- -\frac{1}{w-z_-}.
  $$
  By the definition of the Jacobian of a conformal map:
  $$
  J_{\tau_-^{-1}}(w)= \abs{(\tau_-^{-1})'(w)}^2=\frac{1}{\abs{w-z_-}^4}.
  $$ 
  Hence the sum in~\eqref{e:sumtDl1} takes the  form
  $$
  \sum_{\tDl} \area{\tDl\cap U_-}=\sum_{\tDl}\iint_{\tau_-(\tDl\cap
    U_-)}\frac{1}{\abs{w-z_-}^4}dudv.
  $$
  Since $\{\tau_-(\tDl)\}$ form a perturbed lattice we can estimate
  the right-hand side as an integral of $1/\abs{w-z_-}^4$ over
  $\tau_-(U_-)$ times some universal constant $C$:
  $$
  \sum_{\tDl} \area{\tDl\cap U_-}\ge
  C\iint_{\tau_-(U_-)}\frac{1}{\abs{w-z_-}^4}dudv
  =C\sigma(U_-)=C\pi\lletter^2/2
  $$
  (notice that we used here that $y\in\eta^j(\Delta\cap\RR)$, i.e.,
  the intersection of $U_-$ with $\tDl$ is always non-empty).  Similar
  estimates can be carried over for $U_+$, which yields the inequality
  $$
  \sum_{\tDl} \area{\tDl\cap U}\ge C\pi\lletter^2\ge C\lletter^2.\eqno{\qed}
  $$
  \def\qed{}
\end{proof}


\def\2#1{^{(#1)}}
\def\I#1{I_{#1}}
\def\gm{\gamma}
\def\hyplen(#1){\ell_H(#1)}
\def\gge{>>}
\def\lle{<<}
\section{Uniform twisting on the Julia set}
\label{s:nonlin_Julia}
The main result we prove in this section is the following:
\begin{thm}[{\bf Uniform twisting}]\label{t:nonlin_on_Julia_set}
 There exists a universal constant $C>0$ such that the following statement holds.
  Let $\zeta$ be a  critical commuting pair with an irrational rotation number, and
  let $N$ be such that for all $n>N$ the renormalization $\cR^n \zeta$ extends
  to a $K$-bounded holomorphic commuting pair  with a universal $K$ (as in \S \ref{sec bounds}).
  Fix $n>N$ and let $\cH$ be the above extension.
    Then for every
 $z_0\in{}J(\cH)=\cap_{m\ge0}\cH^{-m}\Delta$ and for all $\lletter\in(0,1]$ there exists
  a pointed area $(U,y)$ such that $\abs{z_0-y}\sim{\lletter}$,
  $\diam{U}\sim{\lletter}$ and
  \begin{equation}
    \label{e:desired_property}
      \sum_{\tilde{\Delta}\subset U}\sg(\tilde{\Delta})\ge Cl^2,
  \end{equation}
  where $\tilde{\Delta}$ is the range of a holomorphic pair, which is
  a universally bounded distortion conformal copy of
  some renormalization of $\cH$, generated by the complete dynamics of
  $\cH$. 

\end{thm}

For $i\geq 0$ set $x_i=f^i(0)$.  Recall the definition of the $m$-th
dynamical partition $\cP_m$ (\ref{dynamical partition}), and for a
point $x\notin\{x_i$, $0\le{}i<q_m+q_{m+1}\}$ denote $\cP_m(x)$ the
element of the partition containing $x$.  For $x=x_i$, let $\cP_m(x)$
be the element of the partition which directly follows $x$ according
to the standard choice of the orientation on the circle.

We will find the following two simple lemmas useful in the proof.
\def\tJ{\tilde J_{m,k}} 
\begin{lem}
  \label{l:f(-i)_Delta_i}
  Let $f:\TT\to\TT$ be a critical circle mapping.  Let $J_m$ be an
  open interval with endpoints $f^{q_m}(0)$ and $f^{q_{m-1}}(0)$. Then
  the intervals $J_m$, $f^{-1}J_m$, \ldots, $f^{-q_{m-1}+1}J_{m}$ are
  all disjoint.  Moreover, there exists $M\in\NN$ such that for all
  $m\geq M$ the following holds.  Let $\tJ$ be the largest closed
  interval containing $f^{-k}J_m$ such that
  $$
  \tJ\subset{}\TT\setminus
  \bigcup_{j\in\{0,\ldots,q_{m-1}-1\}\setminus\{k\}}f^{-j}J_m. 
  $$
  Then $f^{-k}J_m$ lies universally well-inside $\tJ$ and
  $\abs{f^{-k}J_m}\ge\abs{\tJ}/C$ for some universal $C>1$. 
\end{lem}

\begin{proof}
  Let us recall that
  $$
  \cP_{m-2}=
  \left\{
    \begin{array}{ll}
      \I{m-1}^j, & 0\le j<q_{m-2}\\
      \I{m-2}^j, & 0\le j<q_{m-1}\\
    \end{array}
  \right\}
  \text{ and }
  \cP_{m-1}=
  \left\{
    \begin{array}{ll}
      \I{m}^j, & 0\le j<q_{m-1}\\
      \I{m-1}^j, & 0\le j<q_{m}\\
    \end{array}
  \right\}.
  $$
  Notice that $J_m\subset\I{m-1}\cup\I{m-2}$ and
  $f^{-q_{m-1}}J_m\subset{}I_{m-2}$:
\comm{ 
 \begin{center}
    \setlength{\unitlength}{0.00087489in}
    {\renewcommand{\dashlinestretch}{30}
      \begin{picture}(6324,889)(0,-10)
        \put(1204.500,-1040.750){\arc{3220.331}{4.2331}{5.1917}}
        \put(2104.500,-1040.750){\arc{3220.331}{4.2331}{5.1917}}
        \path(12,388)(6312,388)
        \path(462,478)(462,298)
        \path(1362,478)(1362,298)
        \path(1947,478)(1947,298)
        \path(5277,478)(5277,298)
        \put(1362,73){\makebox(0,0)[b]{$0$}}
        \put(462,73){\makebox(0,0)[b]{$x_{q_{m-1}}$}}
        \put(1947,73){\makebox(0,0)[b]{$x_{q_m}$}}
        \put(5277,73){\makebox(0,0)[b]{$x_{q_{m-2}}$}}
        \put(1182,703){\makebox(0,0)[b]{$J_m$}}
        \put(2082,703){\makebox(0,0)[b]{$f^{-q_{m-1}}J_m$}}
      \end{picture}
    }
  \end{center}
}
  Therefore
  \begin{align*}
    f^{-1}J_m&\subset\I{m-2}^{q_{m-1}-1},\\
    \vdots &\\
    f^{-i}J_m&\subset\I{m-2}^{q_{m-1}-i},\\
    \vdots &\\
    f^{-q_{m-1}+1}J_m&\subset\I{m-2}^{q_{m-1}-q_{m-1}+1}=\I{m-2}^1.\\
  \end{align*}
  Hence the intervals $J_m$, $f^{-1}J_m$, \ldots, $f^{-q_{m-1}+1}J_{m}$
  are disjoint. {\it A priori} bounds for $\cP_{m-2}$ imply that the
  adjacent intervals $I'$, $I''\in\cP_{m-2}$ are
  commensurable. {\it A priori} bounds also imply that
  $$
  \frac1{C_0} \abs{\I{m-2}^{q_{m-1}-i}} \le \abs{f^{-i}J_m}\le {C_0}
  \abs{\I{m-2}^{q_{m-1}-i}}
  $$
  The partition  $\cP_{m-2}$ may contain at most
  one interval of the orbit $I_{m-1}^j$ between the interval
  $I_{m-2}^{q_{m-1}-k}\supset{}f^{-k}J_m$ and the
  interval of the orbit $I_{m-2}^j$ which follows it. Therefore the required
  bound for $\tJ$ holds true.
\end{proof}

\begin{lem}
  \label{l:cubic_root_of_the_disc}
  Let $\cH=(\xi,\eta)$ be a universally bounded holomorphic pair. Then there exists a disk
  $D'\subset{}D\setminus\RR$ and a disk $D''\subset\Omega_\cH$,
  centered at the origin, with $\diam{}D'\sim\abs{J}\sim\diam{}D''$
  and $\dist(D',\RR)\sim{}J$, such that $\nu=\xi\circ\eta$ is univalent
  in $D'$ and $D''\subset\nu(D')$.
\end{lem}
\begin{proof}
  Follows from the fact that $\nu$ is cubic up to a universally bounded distortion.
\end{proof}

\noindent
The first step towards proving \thmref{t:nonlin_on_Julia_set} is:
\begin{lem}
  \label{l:nonlin_for_postctitical_set}
  There exists $C>0$ independent of $\zeta$ such that for all $n>N$ as in the 
statement of \thmref{t:nonlin_on_Julia_set}  the following
  statement holds.
    Let $\cH$ be a $K$-bounded holomorphic pair extension of
  $\cR^n\zeta$. For any
  $z_0\in[\xi(0),\eta(0)]$ and for any $\lletter\in(0,1]$ there exists
  a pointed area $(U,y)$ such that $\abs{z_0-y}\sim{\lletter}$,
  $\diam{U}\sim{\lletter}$ and
 the estimate (\ref{e:desired_property}) holds.
\end{lem}

\begin{proof}
  Let $f=\cH|_{[\xi(0),\eta(0)]}$, $\rho(f)=[r_0,r_1,\ldots]$,
  $p_m/q_m=[r_0,\ldots,r_{m-1}]$ and for each $m$ choose ${J}_m$ 
  as in Lemma~\ref{l:f(-i)_Delta_i}. Thus ${J}_m$ is the 
  dynamical interval of the $m$-th pre-renormalization $\cH_m$ of $\cH$.
  Let $\Delta_m\supset{J}_m$ be the range of $\cH_m$. Complex
  bounds~(Theorem~\ref{complex bounds}) imply that
  $\diam\Delta_m\sim\abs{{J}_m}$. For each $0<i\le{}q_{m-1}-1$, let
  $V_m\2i=f^{-i}\Delta_m$. By the Koebe Distortion Theorem we have
  $\diam{V_m\2i}\sim\abs{f^{-i}J_m}$.  Now let us choose minimal $m$ such that
  $$
  \abs{\cP_{m}(z_0)}\le\lletter\le\abs{\cP_{m-1}(z_0)}
  $$
  Consider $r_m=\chi(\cH_m)$, recalling that
  $q_{m+1}=r_mq_m+q_{m-1}$. Consider $m>M$ and take $N$ as in the
  proof of the Lemma~\ref{l:grade_for_almost_parabolic_map}. If
  $r_m\le N$, then 
  $$
  \frac{\abs{\cP_{m}(z_0)}}{\abs{\cP_{m-1}(z_0)}}\ge\frac1{C(N)},\; 
  $$
  uniformly in $z_0$, where $C(N)>1$ is chosen by real {\it a
    priori} bounds and satisfies the additional condition
  $1/C(N)<K_1$, where $K_1>0$ is some universal number to be
  determined in the course of the proof.  Hence for $r_m\le{}N$ we can
  simply choose $k$ so that $f^{-k}J_m$ is the closest to $z_0$. By
  Lemma~\ref{l:f(-i)_Delta_i} and commensurability of $f^{-k}J_m$ and
  $V_m\2k$, this value of $k$ and the domain $(U,y)=(V_m\2k,f^{-k}(0))$
  satisfy the conditions of the Lemma.

  Let us now consider the other case, when $r_m>N$ and the ratio
  $$
  \abs{\cP_{m}(z_0)}/\abs{\cP_{m-1}(z_0)}
  $$
  is smaller when $1/C(N)$. The only complicated case with
  ``possible parabolic cascade'' in this situation is
  \begin{itemize}
  \item[(i)] $\cP_{m-1}(z_0)=\I{m-1}^t$ for some $0\le t<q_m$;
  \item[(ii)] $\cP_{m}(z_0)=\I{m}^s$ for some $0\le s<q_{m+1}$;
  \end{itemize}
  Indeed, it can not happen that $\cP_{m-1}(z_0)=\I{m}^j$,
  $0\le{}j<q_{m-1}$, since then $\cP_{m}(z_0)=\I{m}^j$, and
  $\abs{\cP_{m}(z_0)}/\abs{\cP_{m-1}(z_0)}=1$ and we can
  choose $k$ so that $V_m\2k$ is the closest to $z_0$ as in previous
  paragraphs. 
  
  If $\cP_{m}(z_0)=\I{m+1}^s$, $0\le s<q_m$, then
  $\I{m+1}^s\subset\I{m-1}^s=\cP_{m-1}(z_0)$, i.e., $s=t$ and {\it a priori}
  estimates imply that they are comparable, i.e.,
  $\abs{\cP_{m}(z_0)}/\abs{\cP_{m-1}(z_0)}>K_1$, where $K_1$ is a
  constant whose existence was announced above. Similarly, we can
  choose $k$ so that $V_m\2k$ is the closest to $z_0$, and domain
  $(U,y)=(V_m\2k,f^{-k}(0))$ works.

  Let us sketch the position of intervals in the case of ``possible
  parabolic cascade'':

  \begin{center}
    \setlength{\unitlength}{0.00077489in}
    {\renewcommand{\dashlinestretch}{30}
      \begin{picture}(7674,1035)(0,-10)
        \put(4714.500,1395.000){\arc{1350.750}{1.0475}{2.0941}}
        \put(3612.000,9727.500){\arc{18915.000}{1.2312}{1.9104}}
        \path(12,810)(7662,810)
        \path(462,630)(462,990)
        \path(6762,630)(6762,990)
        \path(4377,855)(4377,765)
        \path(5052,855)(5052,765)
        \path(4782,855)(4692,765)
        \path(4692,855)(4782,765)
        \put(4737,900){\makebox(0,0)[b]{$z_0$}}
        \put(3567,0){\makebox(0,0)[b]{$I_{m-1}^t$}}
        \put(4737,495){\makebox(0,0)[b]{$I_m^s$}}
      \end{picture}
    }
  \end{center}
  
  Applying $f^{-t}$ to $I_{m-1}^t$ we arrive exactly into situation
  described in Lemma~\ref{l:grade_for_almost_parabolic_map}:
  
  \begin{center}
    \setlength{\unitlength}{0.00067366in}
    {\renewcommand{\dashlinestretch}{30}
      \begin{picture}(9024,1174)(0,-10)
        \put(6064.500,1354.000){\arc{1350.750}{1.0475}{2.0941}}
        \put(1542.000,-1499.000){\arc{4986.000}{4.2843}{5.1405}}
        \put(2847.000,-1499.000){\arc{4986.000}{4.2843}{5.1405}}
        \put(7662.000,-311.000){\arc{2340.000}{4.3176}{5.1072}}
        \path(12,769)(9012,769)
        \path(1812,589)(1812,949)
        \path(8112,589)(8112,949)
        \path(5727,814)(5727,724)
        \path(6402,814)(6402,724)
        \path(6132,814)(6042,724)
        \path(6042,814)(6132,724)
        \path(507,589)(507,949)
        \path(2577,589)(2577,949)
        \path(7206,586)(7206,946)
        \path(6942,229)(7167,499)
        \blacken\path(7113.225,387.608)(7167.000,499.000)(7067.131,426.019)(7113.225,434.469)(7113.225,387.608)
        \put(1812,319){\makebox(0,0)[b]{$x_0=0$}}
        \put(2577,274){\makebox(0,0)[b]{$x_{q_{m+1}}$}}
        \put(1542,1039){\makebox(0,0)[b]{$J_{m+1}$}}
        \put(507,319){\makebox(0,0)[b]{$x_{q_m}=\eta_m(0)$}}
        \put(8112,319){\makebox(0,0)[b]{$x_{q_{m-1}}=\xi_m(0)$}}
        \put(2892,1039){\makebox(0,0)[b]{$\eta_m^{-1}J_{m+1}$}}
        \put(6087,859){\makebox(0,0)[b]{$f^{-t}z_0$}}
        \put(6087,424){\makebox(0,0)[b]{$I_m^{s-t}$}}
        \put(7662,949){\makebox(0,0)[b]{$I_m^{q_{m-1}}$}}
        \put(6842,49){\makebox(0,0)[b]{$\eta_m\circ\xi_m(0)=x_{q_{m-1}+q_m}$}}
      \end{picture}
    }
  \end{center}
  
  Notice that $I_m^{s-t}\ne{}I_{m}^{q_{m-1}}$ (otherwise
  $\abs{\cP_{m}(z_0)}/\abs{\cP_{m-1}(z_0)}>K_1$) and
  $I_m^{s-t}\cap{}J_{m+1}=\varnothing$ (by combinatorial reasons).
  
  Let $(\eta_m,\xi_m)=\cR^m(\cH)$.  Consider $J_{m+1}$ and
  corresponding range $\Delta=\Delta_{m+1}$. Obviously,
  $\eta_m^{-1}J_{m+1}\subset{}I_{m-1}$ and $\eta_m^{-2}J_{m+1}$ lies
  well-inside $I_{m-1}$. The same branch $\eta_m^{-1}$ pullbacks
  $\eta_m^{-1}\Delta$ (if $\eta_m^{-1}\Delta\ni0$, we can consider
  $\Delta=\Delta_{m+k}$, such that $\eta_m^{-1}\Delta\not\ni0$, where
  $k$ is universally bounded). By
  Lemma~\ref{l:grade_for_almost_parabolic_map} we can find a pointed
  area $(U',y')$, $\abs{y'-f^{-t}z_0}\sim\lletter$, $\diam
  U'\sim\lletter$, which intersects a large number of images of
  $\eta_m^{-1}\Delta$ under the complete dynamics of $\eta_m$. Moreover,
  $\dist(U',0)\sim\abs{I_{m-1}}$. Hence an application of the Koebe Distortion Theorem
  gives $(U,y)=(f^t(U'),f^ty')$ which satisfies the conditions of the
  theorem.
\end{proof}

By passing to renormalizations of $\zeta$ one can easily derive the following
version of the above lemma:

\begin{prop}
  \label{p:nonlin_for_postctitical_set}
  There exists $C$ such that for all $n$ large enough the following holds.  
  Let $\cH:\Omega\to\Delta$ be a $K$-bounded holomorphic pair extension of
  $\cR^n\zeta$. Then for every
  $z_0\in\Omega\cap\RR$ and for all $\lletter\in(0,1]$ there exists
  a pointed area $(U,y)$ such that $\abs{z_0-y}\sim{\lletter}$,
  $\diam{U}\sim{\lletter}$ and
  (\ref{e:desired_property}) holds.
\end{prop}

\begin{figure}
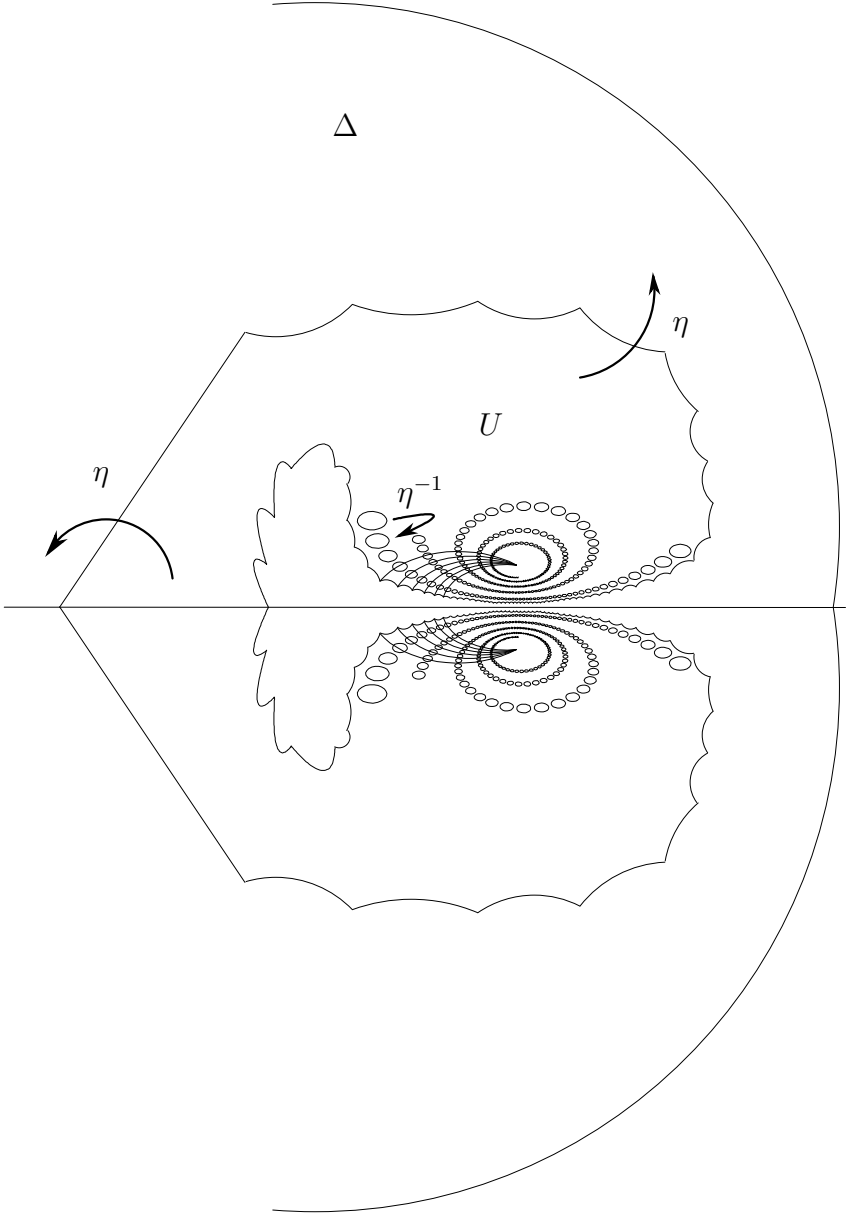

{
\def\phi{\eta}
\parindent=-3cm
\input figpic/wing.eepic

}
\caption{Butterfly wing with parabolic behaviour ($\chi(\cH)>N$)}
\label{f:wing}
\end{figure}

\begin{proof}[Proof of \thmref{t:nonlin_on_Julia_set}]
  The proof is similar to steps II and III of the proof
  of~\cite[Theorem 6.8]{dFdM1} with two essential differences. The
  first difference is that we have to take into account possible
  parabolic cascades and to use
  Lemma~\ref{l:grade_for_almost_parabolic_map} to work with them. The
  second difference is that we should treat specially the cases when
  $\chi(\cH)$ is large and the trajectory $z_k=\cH^kz$ approaches
  the real axis along the wings of the ``butterfly'', not
  near the interval $\Omega_\cH\cap\RR$ (see \figref{pinchfig}); we use
  Lemma~\ref{l:grade_for_almost_parabolic_map} together with
  Lemma~\ref{cut-out} to get the estimates in this case. 
  
  From now on assume $z_0\in J(\cH)\setminus\RR$.  By symmetry it is
  enough to consider the case $z_0\in\HH\cap{}J(\cH)$. Let
  $H=\Delta\setminus\RR$ and consider a vector $v_0$ at $z_0$ with the
  Euclidian norm $\abs{v_0}\sim l$. The sketch of the argument, which goes
  back to McMullen, is as follows: consider the iterations $z_k=\cH^kz_0$,
  $v_k=(\cH^k(z_0))'v_0$ and  wait until the disk of radius $v_k$
  around $z_k$  encloses a commesurable $(U',y')$ disk with
  desired property~\eqref{e:desired_property}. Afterwards we can pull
  $(U',y')$ back to $z_0$ to get the required $(U,y)$. Notice that
  $z_k$ never leaves $J(\cH)$, while the hyperbolic length
  $\hyplen(v_k)$ of $v_k$ increases to infinity (see Theorem 4.11, \cite{dFdM2}):
  $$
  \hyplen(v_0)\le\hyplen(v_1)\le\ldots\le\hyplen(v_k)\to\infty.
  $$
  There are two cases to consider.

  \begin{figure}[htbp]
    \centering
    \setlength{\unitlength}{0.00087489in}
    {\renewcommand{\dashlinestretch}{30}
      \begin{picture}(5500,4795)(0,-10)
        \path(4276,3466)(4186,3556)
        \path(4276,3466)(4186,3556)
        \path(4186,3466)(4276,3556)
        \path(4186,3466)(4276,3556)
        \path(1846,1756)(1756,1846)
        \path(1846,1756)(1756,1846)
        \path(1756,1756)(1846,1846)
        \path(1756,1756)(1846,1846)
        \put(4231,3511){\ellipse{2522}{2522}}
        \put(4231,3511){\ellipse{284}{284}}
        \put(1801,1801){\ellipse{1026}{1026}}
        \put(1801,1801){\ellipse{3586}{3586}}
        \path(181,1711)(2881,1711)
        \dashline{60.000}(2881,1711)(4681,1711)
        \blacken\path(1678.149,1815.408)(1801.000,1801.000)(1699.382,1871.526)(1722.436,1830.727)(1678.149,1815.408)
        \path(1801,1801)(136,2431)
        \blacken\path(258.851,2416.592)(136.000,2431.000)(237.618,2360.474)(214.564,2401.273)(258.851,2416.592)
        \path(1801,1801)(1666,1306)
        \blacken\path(1668.631,1429.665)(1666.000,1306.000)(1726.517,1413.878)(1688.102,1387.040)(1668.631,1429.665)
        \blacken\path(4111.000,3481.000)(4231.000,3511.000)(4111.000,3541.000)(4147.000,3511.000)(4111.000,3481.000)
        \path(4231,3511)(2971,3511)
        \blacken\path(3091.000,3541.000)(2971.000,3511.000)(3091.000,3481.000)(3055.000,3511.000)(3091.000,3541.000)
        \blacken\path(4261.000,3391.000)(4231.000,3511.000)(4201.000,3391.000)(4231.000,3427.000)(4261.000,3391.000)
        \dashline{60.000}(4231,3511)(4231,1711)
        \blacken\path(4201.000,1831.000)(4231.000,1711.000)(4261.000,1831.000)(4231.000,1795.000)(4201.000,1831.000)
        \put(901,2251){\makebox(0,0)[lb]{$R'$}}
        \put(1801,1891){\makebox(0,0)[b]{$z_{k+1}$}}
        \put(1801,1441){\makebox(0,0)[lb]{$v_{k+1}$}}
        \put(4456,3466){\makebox(0,0)[lb]{$z_k$}}
        \put(3601,3601){\makebox(0,0)[lb]{$R$}}
        \put(4321,2656){\makebox(0,0)[lb]{$\Im z_k$}}
      \end{picture}
    }    
    \caption{Case of $\hyplen(v_k)<\eps$ and $\hyplen(v_{k+1})>1/\eps$}
    \label{f:zkzk+1}
  \end{figure}
  
  \emph{Case A}. There exists $k$ such that $\hyplen(v_k)<\eps$ while
  $\hyplen(v_{k+1})>1/\eps$, where $\eps$ is a universal constant to
  be determined in the course of the argument. To be definite, assume
  $z_k$ lies in the domain of $\eta$, so $z_{k+1}=\eta(z_k)$ . Since
  for all $z$ in domain $\Omega$ of $\cH$, the hyperbolic density of
  $\cH$ is commensurable with $1/\Im{z}$, we have
  $\abs{v_k}/\Im{z_k}\sim\hyplen(v_k)<\eps$.  Hence
  $\Im{z_k}\ge{}C\abs{v_k}/\eps\gge\abs{v_k}$. Since all the branching
  points belong to the real axis, $\eta$ is univalent in a disk $D_R(z_k)$
  of radius $R\sim\Im{z_k}$ (notice, that we may need to consider the
  renormalization of $\cH$ instead of $\cH$ to be able to extend
  $\eta$ holomorphically to the disk), cf Figure~\ref{f:zkzk+1}.  By the
  Koebe One-quarter Theorem, $\eta(D_R(z_k))$ contains the disk
  $D_{R'}(z_{k+1})$, where
  $R'=R\abs{v_{k+1}}/4\abs{v_k}\gge\abs{v_{k+1}}$. But
  $\hyplen(v_{k+1})\sim\abs{v_{k+1}}/\Im{z_{k+1}}$ is large
  ($>1/\eps$), hence $\abs{v_{k+1}}\gge\Im{z_{k+1}}$.  Now consider
  two subcases.
  
  \emph{Subcase A$'$}. Assume that $\chi(\cH)\le{}N$, $N$ being a universal
  number to be fixed later. Then the distance of the wings from the
  real axis is greater than some constant $\eps'=\eps'(N)>0$. Hence,
  there exists constant $\eps_1=\eps_1(N)$, such that for any
  $0<\eps<\eps_1$ we can take the point $\zeta$ in the interval
  $\Omega\cap\RR$, closest to $z_{k+1}$, and
  $\dist(\zeta,z_{k+1})\sim\Im{z_{k+1}}$. By
  Proposition~\ref{p:nonlin_for_postctitical_set} there exists a pointed
  domain $(U'',y'')$ with the property~\eqref{e:desired_property},
  such that $\diam{U''}\sim\abs{v_{k+1}}$ and
  $\abs{\zeta-y''}\sim\abs{v_{k+1}}$, and therefore
  $\abs{z_{k+1}-y''}\sim\abs{v_{k+1}}$. If $\eps<\eps_1$ is chosen
  small enough, we have $U''\subset{}D_{R'/2}(z_{k+1})$. Take
  $U'=\eta^{-1}(U'')\subset{}D_R(z_k)$ and $y'=\eta^{-1}y''\in{}U'$.
  By the Koebe Distortion Theorem $\diam(U')\sim\abs{v_k}$, and
  $\abs{z-v_k}\sim\abs{v_k}$. Since $D_R(z_k)\subset{}\HH$, we can
  pull back $(U',y')$ using the univalent inverse branch of $\cH^{-k}$
  mapping $z_k$ to $z_0$ to get $(U,y)=(\cH^{-k}U',\cH^{-k}y')$, which
  by the Koebe Theorem satisfies~\eqref{e:desired_property} with a
  finitely distorted constant $C$.

  \begin{figure}
    \begin{center}
      \setlength{\unitlength}{0.00077489in}
{\renewcommand{\dashlinestretch}{30}
\begin{picture}(7224,2212)(0,-10)
\put(2541.000,495.000){\arc{900.000}{3.7851}{5.6397}}
\put(1800.107,453.214){\arc{984.459}{3.9524}{5.5972}}
\put(1080.107,498.214){\arc{984.459}{3.9524}{5.5972}}
\put(346.547,624.375){\arc{913.851}{3.9194}{5.7541}}
\put(3261.000,495.000){\arc{900.000}{3.7851}{5.6397}}
\put(4686.000,495.000){\arc{900.000}{3.7851}{5.6397}}
\put(5426.893,453.214){\arc{984.459}{3.8276}{5.4724}}
\put(6146.893,498.214){\arc{984.459}{3.8276}{5.4724}}
\put(6880.453,624.375){\arc{913.852}{3.6707}{5.5054}}
\put(3966.000,495.000){\arc{900.000}{3.7851}{5.6397}}
\put(2541.000,940.000){\arc{900.000}{0.6435}{2.4981}}
\put(1800.107,981.786){\arc{984.459}{0.6860}{2.3308}}
\put(1080.107,936.786){\arc{984.459}{0.6860}{2.3308}}
\put(346.547,810.625){\arc{913.852}{0.5291}{2.3638}}
\put(3261.000,940.000){\arc{900.000}{0.6435}{2.4981}}
\put(4686.000,940.000){\arc{900.000}{0.6435}{2.4981}}
\put(5426.893,981.786){\arc{984.459}{0.8108}{2.4556}}
\put(6146.893,936.786){\arc{984.459}{0.8108}{2.4556}}
\put(6880.453,810.625){\arc{913.851}{0.7778}{2.6125}}
\put(3966.000,940.000){\arc{900.000}{0.6435}{2.4981}}
\put(3343,1980){\ellipse{360}{180}}
\put(561,2100){\ellipse{360}{180}}
\put(2643,1985){\ellipse{360}{180}}
\put(1933,2010){\ellipse{360}{180}}
\put(1236,2048){\ellipse{360}{180}}
\put(4048,1980){\ellipse{360}{180}}
\put(6830,2100){\ellipse{360}{180}}
\put(4748,1985){\ellipse{360}{180}}
\put(5458,2010){\ellipse{360}{180}}
\put(6155,2048){\ellipse{360}{180}}
\put(3261,1350){\ellipse{360}{180}}
\put(479,1470){\ellipse{360}{180}}
\put(1851,1380){\ellipse{360}{180}}
\put(1154,1418){\ellipse{360}{180}}
\put(3966,1350){\ellipse{360}{180}}
\put(6748,1470){\ellipse{360}{180}}
\put(4666,1355){\ellipse{360}{180}}
\put(5376,1380){\ellipse{360}{180}}
\put(6073,1418){\ellipse{360}{180}}
\put(2561,1355){\ellipse{360}{180}}
\put(2902,796){\ellipse{1080}{1080}}
\path(12,720)(7212,720)
\path(2857,751)(2947,841)
\path(2857,841)(2947,751)
\path(2907,791)(2907,1336)
\blacken\path(2937.000,1216.000)(2907.000,1336.000)(2877.000,1216.000)(2907.000,1252.000)(2937.000,1216.000)
\put(2972,1026){\makebox(0,0)[lb]{$v_{k+1}$}}
\put(3032,731){\makebox(0,0)[lb]{$z_{k+1}$}}
\put(2422,1476){\makebox(0,0)[lb]{$\tilde\Delta$}}
\put(2302,901){\makebox(0,0)[rb]{$\gamma^i_m$}}
\end{picture}
}

      ``Degenerate case''

      \bigskip
      \bigskip
      \bigskip

      \setlength{\unitlength}{0.00037489in}
{\renewcommand{\dashlinestretch}{30}
\begin{picture}(11570,8121)(0,-10)
\put(4688.000,2883.000){\arc{900.000}{3.7851}{5.6397}}
\put(3947.107,2841.214){\arc{984.459}{3.9524}{5.5972}}
\put(3227.107,2886.214){\arc{984.459}{3.9524}{5.5972}}
\put(2493.547,3012.375){\arc{913.851}{3.9194}{5.7541}}
\put(5408.000,2883.000){\arc{900.000}{3.7851}{5.6397}}
\put(1773.547,3102.375){\arc{913.851}{3.9194}{5.7541}}
\put(1055.272,3206.182){\arc{897.208}{3.8947}{5.7787}}
\put(324.875,3350.500){\arc{869.289}{3.8955}{5.9000}}
\put(6833.000,2883.000){\arc{900.000}{3.7851}{5.6397}}
\put(7573.893,2841.214){\arc{984.459}{3.8276}{5.4724}}
\put(8293.893,2886.214){\arc{984.459}{3.8276}{5.4724}}
\put(9027.453,3012.375){\arc{913.852}{3.6707}{5.5054}}
\put(6113.000,2883.000){\arc{900.000}{3.7851}{5.6397}}
\put(9747.453,3102.375){\arc{913.852}{3.6707}{5.5054}}
\put(10465.728,3206.182){\arc{897.205}{3.6460}{5.5300}}
\put(11196.125,3350.500){\arc{869.289}{3.5248}{5.5293}}
\put(4695.000,2512.000){\arc{900.000}{0.6435}{2.4981}}
\put(3954.107,2553.786){\arc{984.459}{0.6860}{2.3308}}
\put(3234.107,2508.786){\arc{984.459}{0.6860}{2.3308}}
\put(2500.547,2382.625){\arc{913.852}{0.5291}{2.3638}}
\put(5415.000,2512.000){\arc{900.000}{0.6435}{2.4981}}
\put(1780.547,2292.625){\arc{913.852}{0.5291}{2.3638}}
\put(1062.272,2188.818){\arc{897.205}{0.5044}{2.3884}}
\put(331.875,2044.500){\arc{869.289}{0.3832}{2.3877}}
\put(6840.000,2512.000){\arc{900.000}{0.6435}{2.4981}}
\put(7580.893,2553.786){\arc{984.459}{0.8108}{2.4556}}
\put(8300.893,2508.786){\arc{984.459}{0.8108}{2.4556}}
\put(9034.453,2382.625){\arc{913.851}{0.7778}{2.6125}}
\put(6120.000,2512.000){\arc{900.000}{0.6435}{2.4981}}
\put(9754.453,2292.625){\arc{913.851}{0.7778}{2.6125}}
\put(10472.728,2188.818){\arc{897.208}{0.7532}{2.6372}}
\put(11203.125,2044.500){\arc{869.289}{0.7539}{2.7584}}
\put(5490,4368){\ellipse{360}{180}}
\put(563,4773){\ellipse{360}{180}}
\put(2708,4488){\ellipse{360}{180}}
\put(1995,4563){\ellipse{360}{180}}
\put(1298,4661){\ellipse{360}{180}}
\put(4790,4373){\ellipse{360}{180}}
\put(4080,4398){\ellipse{360}{180}}
\put(3383,4436){\ellipse{360}{180}}
\put(6195,4368){\ellipse{360}{180}}
\put(11122,4773){\ellipse{360}{180}}
\put(8977,4488){\ellipse{360}{180}}
\put(9690,4563){\ellipse{360}{180}}
\put(10387,4661){\ellipse{360}{180}}
\put(6895,4373){\ellipse{360}{180}}
\put(7605,4398){\ellipse{360}{180}}
\put(8302,4436){\ellipse{360}{180}}
\put(5408,3738){\ellipse{360}{180}}
\put(481,4143){\ellipse{360}{180}}
\put(2626,3858){\ellipse{360}{180}}
\put(1913,3933){\ellipse{360}{180}}
\put(1216,4031){\ellipse{360}{180}}
\put(4708,3743){\ellipse{360}{180}}
\put(3998,3768){\ellipse{360}{180}}
\put(3301,3806){\ellipse{360}{180}}
\put(6113,3738){\ellipse{360}{180}}
\put(11040,4143){\ellipse{360}{180}}
\put(8895,3858){\ellipse{360}{180}}
\put(9608,3933){\ellipse{360}{180}}
\put(10305,4031){\ellipse{360}{180}}
\put(6813,3743){\ellipse{360}{180}}
\put(7523,3768){\ellipse{360}{180}}
\put(8220,3806){\ellipse{360}{180}}
\put(4313,4053){\ellipse{8090}{8090}}
\path(1208,2748)(11558,2748)
\path(4268,4143)(4403,4008)
\path(4268,4008)(4403,4143)
\path(4313,4053)(7103,1173)
\blacken\thicklines
\path(6892.915,1303.630)(7103.000,1173.000)(6979.104,1387.125)(6986.106,1293.664)(6892.915,1303.630)
\put(4493,4008){\makebox(0,0)[lb]{$z_{k+1}$}}
\put(6068,2388){\makebox(0,0)[lb]{$v_{k+1}$}}
\end{picture}
}

      ``Non-degenerate case''
    \end{center}

    \caption{$z_{k+1}$ and parabolic cascade in wings}
    \label{f:zk1parab}
  \end{figure}
  
  \emph{Subcase A$''$}. $\chi(\cH)>N$. This case is illustrated by
  Figure~\ref{f:wing}. Lemma~\ref{cut-out} guarantees that the mapping
  $\eta$ is close to a parabolic one in a disk $\hat D$, which is commensurable
  with $I_\eta$. At the same time since $\eta$  also has a branching
  point at the endpoint of the interval $\Omega\cap\RR$, we can bring a
  copy of holomorphic pair, say $\cR\cH$, commesurable with the length
  of $\gamma_m^i$ inside $\hat{D}$ and  use the method of the proof of
  Lemma~\ref{l:grade_for_almost_parabolic_map} to spread it around
  taking images and pre-images with $\eta$. Obviously, if $z_{k+1}$ is
  outside $\hat{D}$ it can be treated as in Subcase $A'$.
  
  Assume that $z_{k+1}$ is inside $\hat{D}\cap{}U$. In order to apply
  estimates similar to the one used in the proof of
  Proposition~\ref{p:nonlin_for_postctitical_set} we have to guarantee
  that for $\tilde\Delta$ closest to $z_{k+1}$ the $\diam\tilde\Delta$
  is smaller than some constant times $\abs{v_{k+1}}$.

  Recall that $\abs{v_{k+1}}\gge\Im{z_{k+1}}$, hence we fall into one
  of the situations, sketched on Figure~\ref{f:zk1parab}. In
  ``nondegenerate case'' $z_{k+1}$ lies inside the lattice of disks,
  which together with Lemma~\ref{l:grade_for_almost_parabolic_map}
  implies that disk $D_{R'}(z_{k+1})$ encloses pointed region
  $(U'',y'')$ with the property~\eqref{e:desired_property} and further
  argument is the same as in Subcase $A'$.
  
  The degenerate situation, sketched in Figure~\ref{f:zk1parab} has,
  in fact, bounded geometry. Indeed, if $z_{k+1}$ lies outside the
  lattice of disks, consider disk $\tilde\Delta$ and a piece of
  boundary $\gm_m^i$, both closest to $z_{k+1}$ (cf
  Figure~\ref{f:wing}). From estimates, given in~\cite{Ya1}, and
  boundedness of derivatives of $\cE_\eta^{-1}$ we obtain
  $\dist(\gm_m^i,\RR)\sim\dist(\gm_m^i,\tilde\Delta)\sim\diam\tilde\Delta$
  universally.  Hence $\Im{}z_{k+1}\ge K_0\dist(\gm_m^i,\RR)$, and
  hence $\abs{v_{k+1}}\ge{}K_0\Im{}z_{k+1}$ and one can take
  $(\tilde\Delta,y'')$, $y''\in\Delta$, as a pointed domain
  $(U'',y'')$ with property~\eqref{e:desired_property}, such that
  $\diam{U''}\sim\abs{v_{k+1}}$ and
  $\abs{z_{k+1}-y''}\sim\abs{v_{k+1}}$. Now the same procedure as in
  Subcase $A'$ can be followed (possibly with a further correction of
  $\eps$).

  \emph{Case B}. $\hyplen(v_k)\sim1$ for some $k$. Remember that
  $\Im{z_{k+1}}\sim\abs{v_k}$.
  
  \emph{Subcase $B'$}. $\chi(\cH)>N$ and $z_k\in\hat{D}$. If $z_k$ is
  above the repelling fixed point $z_+$ (i.e. $\Im{}z_k\ge\Im{}z_+$),
  then by estimates similar to the one used in the proof of
  Lemma~\ref{l:nonlin_for_postctitical_set} there exists a domain
  $(U',y')$ with property~\eqref{e:desired_property} such that
  $\diam{U'}\sim\abs{v_{k+1}}$ and
  $\abs{z_{k+1}-y'}\sim\abs{v_{k+1}}$. This $(U',y')$ is now easily
  pulled back by $\cH^{-k}$. If $z_k$ is below the fixed point $z_+$,
  then there are two subcases depending on the $l=\dist(z_k,\gm_m^i)$,
  where $\gm_m^i$ is the closest piece of $\hat{D}\cap\Omega$ to
  $z_k$.
  
  These subcases are similar to the ``non-degenerate'' and ``degenerate''
  situations, considered in Subcase $A''$ with only the difference that in
  the ``degenerate'' situation one can not hope for a large disk
  $D_{R'}(z_{k+1})$, enclosing many disks from the parabolic lattice
  (because $R'$ is comparable with $\abs{v_{k+1}}$ here). Hence, in
  the ``degenerate'' situation we iterate $z_{k+1}$ further, until such a
  disk would emerge. The exact description of dynamics, given by
  Lemma~\ref{cut-out} allows us to do so. Next, we describe this
  procedure more formally.
  
  If $l>K_1\abs{v_k}$ for some $K_1$, then we can use the same
  estimates as in Lemma~\ref{l:nonlin_for_postctitical_set} and the
  argument is similar to the one used in subcase $A''$. Otherwise
  $z_k$ becomes quite close to $\gm_m^i$ and further application of
  $\eta$ moves $z_k$ along the parabolic cascade, increasing the
  hyperbolic length of $\hyplen(v_k)$ (see Lemma~\ref{cut-out} for
  description of dynamics).

  Eventually, we either fall into Subcase $A''$ (with the large
  hyperbolic length and the possible pull-back of the large number of
  copies of $\tilde\Delta$), or $\hyplen(v_l)$ remains comparable with
  1 (which also allows to pull back a disk with copies of
  $\tilde\Delta$ at the moment $l$ when $z_l\in\hat{D}$,
  $z_{l+1}\notin\hat{D}$).

  \emph{Subcase $B''$}. Situations (1) $\chi(\cH)\le{}N$ and (2)
  $\chi(\cH)>{}N$ and $z_k\notin\hat{D}$ can be considered simultaneously.

  Let $\zeta$ be the point in $\Omega\cap\RR$, closest to
  $z_{k}$. Then by Proposition~\ref{p:nonlin_for_postctitical_set} there
  exists a pointed domain $(U'',y'')$, $y''\in\RR$, such that
  $\abs{\zeta-y''}\sim\Im{}z_k$, $\diam{U''}\sim\Im{}z_k$, satisfying the
  condition~\eqref{e:desired_property}. It follows from the proof of
  Lemma~\ref{l:nonlin_for_postctitical_set} that one of the following
  conditions holds for some universal $K_2$:

  i) There exists $\tilde\Delta\subset{}U''$ such that
  $\diam\tilde\Delta\ge{}K_2\diam{U''}$,
  $\tilde\Delta\cap\RR\ne\varnothing$. 

  ii) There exists a disk $U'$ around $y'$, $U'\subset U''$ , with the
  property~\eqref{e:desired_property}, such that $\diam{U'}\ge{}K_2\diam{U''}$
  and $\dist(U',\RR)\sim\diam{}U''$.
  
  In case i) we can pull-back $\tilde\Delta$ to origin by $\cH^{-k}$
  for some $k$ and to produce a disk $D'$, containing a preimage of a
  holomorphic pair, commensurable with $D'$, using
  Lemma~\ref{l:cubic_root_of_the_disc}. Now we choose
  $(U',y)=(\cH^{k}D',\cH^{k}d)$, where $d$ is the center of $D'$.
  Notice that by the Koebe Theorem $U'$ has the
  property~\eqref{e:desired_property}.

  In case ii) we already have $U'$
  satisfying~\eqref{e:desired_property}.

  The remaining step is to pull-back $(U',y')$ to $z_0$, which can be
  carried over as in~\cite[Theorem 6.8, step IIIb]{dFdM1}. Indeed, if
  $\gm$ is a geodesic arc in $H$, connecting $z_k$ to $U'$, then
  $\diam(\gm\cup{}U')\sim\Im{}z_k$. Hence there exist $r_k\sim\Im{}z_k$
  and $\zeta_1=z_k$, $\zeta_2$, \ldots, $z_s\in\gm\cup{}U'$, such that
  $$
  \gm\cup{}U\subset{\Nu}=\cup_{j=1}^s D_{r_k/2}(\zeta_j),
  $$
  and such that $D_{r_k}(\zeta_j)\subset\Delta\setminus\RR$ for all
  $j$. Inverse branch $\Psi=\cH^{-k}$, taking $z_k$ to $z_0$ is
  defined in each $D_{r_k}(\zeta_j)$. Application of the Koebe Theorem to
  each $D_{r_k}(\zeta_j)$ yields
  $\abs{\Psi'(z)}\sim\abs{\Psi'(z_k)}\sim\abs{v_0}/\abs{v_k}$. Hence
  $(U,y)=(\Psi(U'),\Psi(y'))$ satisfy $\diam{}U\sim\abs{v_0}$,
  $\dist(y,z_0)\sim\abs{v_0}$.
\end{proof}

\begin{proof}[Proof of Theorem~\ref{uniform twisting}]
 In view of \propref{dFdM-prop}, the abundance of copies of holomorphic commuting pairs given by the 
estimate~(\ref{e:desired_property}) provides a uniform and universal
  bound from below on nonlinearity defined by~(\ref{nonlinearity}).
\end{proof}

\begin{proof}[Proof of \thmref{main theorem}]
The statement follows from the Dynamical Inflexibility \thmref{thm-inflex} together with
\thmref{twisting bounded type} (I) and 
\thmref{uniform twisting}.
\end{proof}


\section{Hyperbolicity of the renormalization horseshoe}
In this section we will apply our Rigidity Theorem to obtain a new proof of the hyperbolicity of the global
horseshoe for the cylinder renormalization operator  $\cren$ \cite{Ya3,Ya4}. The second author has established this
result in \cite{Ya4} using some infinite-dimensional quasiconformal deformation spaces arguments. Now we will be able to
give a simpler argument, following along the same lines as the proof of hyperbolicity of the periodic orbits of 
$\cren$ given in \cite{Ya3}. We will use \cite{Ya3} as the general reference for this section, however, for the sake
completeness, we will briefly recall the basic definitions. 

\medskip

\noindent
{\bf Definition of $\cren$.}
The main point of \cite{Ya3} was to replace the
renormalization operator $\cR$ acting on the space of commuting pairs with an analytic operator $\cren$
defined on a complex-analytic Banach manifold. To define $\cren$, we will need a few preliminaries.
Firstly, let us denote $\pi:\CC\to\CC/\ZZ$ the natural projection. For an equatorial topological annulus
$U\subset\CC/\ZZ$ denote ${\aaa A}_U$ the space of bounded analytic maps $\phi:U\to\CC/\ZZ$,
such that $\phi(\TT)$ is homotopic to $\TT$, equipped with the uniform metric.
To turn this space into a Banach manifold, consider the  Banach space $\tl{\aaa A}_U$ of bounded
analytic $1$-periodic functions from $\pi^{-1}(U)\to\CC$ with the sup norm, and use the local
homeomorphism $\tl{\aaa A}_U\to{\aaa A}_U$ given by
$$\psi\mapsto\pi\circ(\psi+\text{Id})\circ\pi^{-1}$$
to define the atlas on ${\aaa A}_U$. We denote ${\aaa C}_U$ the codimension two submanifold
of ${\aaa A}_U$ consisting of maps with a cubic critical point at the origin.
Let $\cur\subset \cu$ be the real Banach manifold consisting of the critical circle maps
in ${\aaa C}_U$ (the {\it real slice} of $\cu$).
As in \cite{Ya3}, a  tangent space to $\cu$ will be 
 naturally identified with a  Banach subspace ${\aaa B}_U\subset\tl{\aaa A}_U$,
 ${\aaa B}_U^\RR$ will again denote the real slice.

Given a critical cylinder map $f\in\cu$ let us say that it is {\it cylinder renormalizable}, or simply
{\it renormalizable},
if there exists $k>1$ and an equatorial annulus $V\subset \CC/\ZZ$ such that following holds:
\begin{itemize}
\item there exist repelling periodic points $p_1$, $p_2$ of $f$ in $U$ with periods $k$
and a simple arc $l$ connecting them such that $f^k(l)$ is a simple arc, and  $f^k(l)\cap l=\{p_1,p_2\}$;
\item the iterate $f^k$  is defined and univalent in the domain $C_f$  bounded by $l$ and $f^k(l)$,
the corresponding inverse branch $f^{-k}|_{f^k(C_f)}$ univalently extends to $C_f$;
and the quotient of $\overline{C_f\cup f^k(C_f)}\sm\{p_1,p_2\}$ by the action of $f^k$ is a 
Riemann surface conformally isomorphic to the cylinder $\CC/\ZZ$ 
(we will  call a domain $C_f$ with these properties
a {\it fundamental crescent of $f^k$});
\item for a point $z\in \bar C_f$ with $\{f^j(z)\}_{j\in\NN}\cap \bar C_f\ne \emptyset$,
set $R_{C_f}(z)=f^{n(z)}(z)$ where $n(z)\in\NN$ is the smallest value for which $f^{n(z)}(z)\in\bar C_f$.
We further require that 
there exists a point $c$ in the domain of $R_{C_f}$ such that $f^m(c)=0$ for some $m<n(c)$;
and if we denote $\hat f$ the projection of $R_{C_f}$ to $\CC/\ZZ$ with $c\mapsto 0$, then
$\hat f\in \cv$.
\end{itemize}
We will say that the new critical circle map $\hat f$ is a {\it cylinder renormalization} of $f$
with period $k$.

\begin{prop}[{\bf Properties of cylinder renormalization}, \cite{Ya3}]
\label{ren open set}
The following statements hold:
\begin{itemize}
\item The cylinder renormalization does not depend on  the choice of the fundamental crescent.
\item If $f\in \cu$ is cylinder renormalizable with period $k$, then there exists an open neighborhood
$\cW(f)\subset\cu$ such that every $g\in\cW(f)$ is also cylinder renormalizable, with the same period, and
the fundamental crescent $C_g$ can be chosen to move continuously with $g$.
\item Moreover, the cylinder renormalization is an analytic operator $\cW(f)\to\cv$.
\end{itemize}
\end{prop}

\noindent
The connection with the renormalization of commuting pairs is established as follows:
\begin{prop}[{\bf Cylinders of commuting pairs in $\cE$ } \cite{Ya3}]
\label{cyl-eps}
Let $\zeta=(\eta,\xi)$ be a commuting pair in the Epstein class with $\chi(\zeta)\neq\infty$. Then the map $\eta$
has a real-symmetric fundamental crescent $C_\eta$, with 
$$C_\eta\cup\overline{\eta(C_\eta)}/\eta\simeq\CC/\ZZ,$$
 whose first return map projects to an analytic critical circle map $f_\zeta$
independent of the choice of $C_\eta$, whose rotation number $\rho(f_\zeta)=\rho(\cR\zeta)$. Further, there exists
$U=U(s)$ such that if $\zeta\in \cE_s$, then $f_\zeta\in\cur$.
\end{prop}

\begin{prop}[\cite{Ya3}]
\label{dist}
If $\zeta_1$, $\zeta_2\in\cE$, and $\chi(\zeta_1)=\chi(\zeta_2)=r\neq\infty$, then the equality
$$f_{\zeta_1}\equiv f_{\zeta_2}$$
is equivalent to the existence of a conformal conjugacy between $\zeta_1$ and $\zeta_2$ whose domain 
contains a fundamental crescent (we will write $\zeta_1\ceq\zeta_2$ in this case).
 Moreover, if $\zeta_1,\zeta_2\in\cE_s$, then
$$\dist_{C^0}(\zeta_1,\zeta_2)\geq c(s,r)\dist_{C^0}(f_{\zeta_1},f_{\zeta_2}).$$
\end{prop}

From now on, let us fix $s>0$ as in \lemref{real-bounds}, and set $U=U(s)$ from \propref{cyl-eps}.
Finally, we have:
\begin{prop}
\label{choice N}
Then there exists $N\in\NN$ such that the following holds.
Suppose $\zeta\in\cE_s$ is at least $N$-times renormalizable. Then 
the renormalization $\cR^{N}\zeta\in\cE_s$. Further, the map $f_\zeta$ is cylinder renormalizable
in such a way that its renormalization is 
$\hat f=f_{\cR^{N}\zeta}$, which is in $\cv$ with $V\Supset U$.
\end{prop}

\begin{defn}\label{defn of cyl ren}
For the remainder of the paper fix the value of $N$ as above, and 
set $\mfld=\cur$.
We will call $\hat f$ as in \propref{choice N}  {\em the}  cylinder renormalization of $f_\zeta$, and write 
$$\hat f\equiv\cren f_\zeta.$$
By \propref{ren open set},
for every pair $\zeta$ as above, the tranformation $f_\zeta\mapsto \cren f_\zeta$ extends
to an open neighborhood $Y\subset \cu$ as an analytic operator $Y\to\cu$.
We shall call this operator {\it the cylinder renormalization operator.}
\end{defn}

\noindent
\begin{prop}
\label{RN projects}
If $\zeta_1$, $\zeta_2$ are at least $N$-times renormalizable elements of $\cE_s$,
and $\zeta_1\ceq\zeta_2$, then $\cR^N\zeta_1\ceq\cR^N\zeta_2$. Hence
the action of $\cR^N$ is well-defined on the quotient space $\cE_s/\ceq$, and denoting $\iota:[\zeta]_{\ceq}\to\cur$,
we have
$$\iota\circ \cR^N=\cren\circ\iota.$$ 
\end{prop}

\medskip
\noindent
{\bf Hyperbolicity of the renormalization horseshoe.}
Now let $\hat f\in\mfld$ be a point in the renormalization horseshoe. That is, $\hat f=f_\zeta$ for a commuting pair $\zeta$ in the
horseshoe $\cI$ of $\cR$ (see \cite{Ya2}). Set $\rho=\rho(\hat f)$ and define 
$$\cD_\rho=\{ f\in\mfld,\text{ such that }\rho(f)=\rho\}.$$
We have:
\begin{thm}[Theorem 8.2, \cite{Ya3}]
There exists an open neighborhood $W\subset \mfld$ of $\hat f$ such that $\cD_\rho\cap W$ is a smooth submanifold of
$\mfld$ of codimension $1$.
\end{thm}

\noindent
Let the hyperplane  $T\equiv T_{\hat f}(\cD_\rho\cap W)$ be the 
tangent space to this codimension one submanifold 
at $\hat f$, and let $\cL$ be the differential of $\cren$ at $\hat f$:
$$\cL=D_{\hat f}\cren^m:{\aaa B}_U^\RR\to {\aaa B}_U^\RR,\; \cL:T\to T$$
Recall that a continuous linear operator on a Banach space is called {\it compact} if it maps the 
closed unit ball of the space onto
a compact set. This condition is equivalent to the image of every closed bounded set being compact.
\begin{prop}[Prop. 9.1, \cite{Ya3}]
\label{L compact}
The operator $\cL=D_{\hat f}\cren^m:{\aaa B}_U^\RR\to {\aaa B}_U^\RR$ is compact.
\end{prop}

\noindent
\begin{prop}
\label{prop-cont}
There exists  $a<1$ independent of $\hat f$ such that 
the operator $\cL|_T$ is a contraction by $a$, and moreover, 
the spectral radius $R_{\text{sp}}$ of the operator $\cL|_T$ is at most $a$.
\end{prop}
\begin{pf}
The spectral theory of compact operators implies that there exist finitely many eigenvalues $\lambda_i$ of $\cL|_T$ with
$|\lambda_i|=R_{\text{sp}}(\cL|_T)$.
The space $T$ breaks into a direct sum $E^1\oplus E^2$, where the latter is a finite-dimensional subspace spanned by the
generalized eigenvectors of $\lambda_i$, and the action of $\cren|_{\cD_\rho}$ near $\hat f$ is dominated by 
the projection of $\cL|_T$ to $E^2$. 
The Rigidity Theorem implies that $\cren$ is a universally geometric contraction in $W\cap \cD_\rho$ (\corref{ren conv}). This
implies that $R_{\text{sp}}(\cL|_T)<a<1.$
\end{pf}

Define a cone  $\cC\in{\aaa B}_U^\RR$  as follows:
$$\cC_f=\{v\in {\aaa B}_U^\RR\text{ such that }\inf_{x\in\RR}v(x)>0\}.$$

\noindent
\begin{lem}
\label{invt cone field}

\noindent
\begin{itemize}
\item[(I)] The cone  $\cC$ is renormalization-invariant: 
$\cL:\cC\to \cC$,
\item[(II)] Moreover, there exists $\alpha>0$ and $k\in\NN$ independent of $\hat f$ such that for any vector field  $v\in \cC$
$$\inf_{x\in\RR}D_{\hat f}\cren^{mk}(v(x))>(1+\alpha)\inf_{x\in\RR}v(x)$$
\item[(III)] Finally, 
there exists $\ell\in \NN$ such that 
if $v\in {\aaa B}_U$ belongs to the closure of $\cC$, and $v\ne 0$, then 
$\cL^\ell(v)\in\cC$.
\end{itemize}
\end{lem}

\begin{pf}
Let $\tl f=\pi^{-1}\circ \hat f\circ \pi$. 
As follows from an elementary computation, if $\tl f_t(x)=\tl f(x)+tv(x)+o(t)$,
then 
\begin{equation}
\label{e1}
\tl f_t^2(x)=\tl f^2(x)+tv_2(x)+o(t)=\tl f^2(x)+t(\tl f'(\tl f(x))v(x)+v(\tl f(x)))+o(t),
\end{equation}
and more generally, if we write
\begin{equation}
\label{e2}
\tl f_t^n(x)=\tl f^n(x)+tv_n(x)+o(t),\text{ then }v_n(x)=\tl f'(\tl f^{n-1}(x))v_{n-1}(x)+v(\tl f^{n-1}(x))
\end{equation}
Let us denote $C_k\Subset U$ the fundamental crescent of $\hat f$ corresponding to the 
$k$-th cylinder renormalization $\cren^k \hat f$, and $J_k=C_k\cap\TT$. Let $\phi_k(x)$ be the
corresponding
uniformizing coordinate $C_k\to\CC/\ZZ$, and $\Phi_k:C_n\to\CC$ its lift. In the 
local chart given by $\phi_k^{-1}$, the $k$-th cylinder renormalization $\cren^k \hat f$
is represented by the iterate $g_k=\hat f^{q_{mk}}$ for some $m\in\NN$.

To prove the first claim, observe that by (\ref{e2}),
\begin{equation}
\label{e3}
\inf_{x\in\RR}v_n(x)\geq \inf_{x\in\RR}v(x)>0.
\end{equation} 
The image
\begin{equation}
\label{e4}
D\cren v=[(\Phi_k)'\circ g_k\cdot (v_{q_{mk}}|_{J_k})]\circ (\Phi_k)^{-1}
\end{equation}
Since $(\Phi_k)'>0$, (I) follows.

To prove (II), note that by the real {\it a priori} bounds, the sizes of the intervals $J_k$
decrease geometrically with $k$. On the other hand, applying the Koebe Distortion Theorem
to the conformal extension of $\Phi_k$ to $C_k$ and its two neighboring iterates,
we see that the distortion of $\Phi_k$ on the interval $g_k(J_k)$ is bounded uniformly in $k$ nd
$f$.
Hence, there exists $\alpha>0$ and $l\in\NN$  independent of $f$ such that for all $k\geq l$ 
$$(\Phi_k)'|_{g_k(J_k)}>1+\alpha,$$
and (II) follows from (\ref{e4}).

Finally, to see (III), observe that if $v(x)\not\equiv 0$, then there exist $\ell$, $n\leq q_{m\ell}$
such that $v(x)>0$ for every $x\in f^n(J_\ell)$. The claim now follows from (\ref{e2}).
\end{pf}

\noindent
We conclude:
\begin{prop}
\label{exp dir}
There exists $\ell\in\NN$ such that the invariant horseshoe of the operator $\cren^\ell$ 
is hyperbolic with one-dimensional unstable direction.
\end{prop}
\begin{pf}
Let $k$, $m$, $\alpha$ be as in the previous lemma, and set $\ell=mk$.
In view of \propref{prop-cont} we only need to construct the one-dimensional expanding direction of
$\cren^\ell$. For a periodic point $\hat f$ of period $p$ of this operator, the previous lemma implies that
$R_\text{sp}(D_{\hat f}\cren^{\ell p})>(1+\alpha)^p.$ Moreover, since the spectrum of a compact operator
is discrete, and every non-zero element of the spectrum is an eigenvalue, the corresponding eigendirections
vary continuously with $\hat f$. The statement easily follows.
\end{pf}

\section{Concluding remarks}

\noindent
Let us conclude by outlining some of the recent developments, and open
problems which naturally relate to the results of this paper. Firstly,
de~Faria and de~Melo \cite{dFdM1} have introduced a class of rotation
numbers $\cF\cM\subsetneq \TT\setminus\QQ$ specified by the following
conditions: $\rho=[r_0,r_1,\ldots]\in\cF\cM$ if
\begin{itemize}
\item $\lim\sup \frac{1}{n}\sum_{j=0}^{n-1}\log r_j<\infty$,
\item $\lim\frac{1}{n}\log r_n=0$,
\item
  $\frac{1}{n}\sum_{j=k}^{k+n-1}\log{}r_j\leq{}\omega_\rho\left(\frac{n}{k}\right)$,
  for $0< n\leq k$, where $\omega_\rho(t)$ is a positive function
  defined for $t>0$ and such that $t\omega_\rho(t)\to 0$ as $t\to 0$.
\end{itemize}
\noindent
The class $\cF\cM$ has full Lebesgue measure on the circle, and
contains all the rotation numbers of bounded type.  de~Faria and
de~Melo proved the following:

\medskip

\noindent
{\bf Theorem \cite{dFdM1}} {\sl Let $f$ and $g$ be two critical circle
  mappings of class $C^3$ with the same irrational rotation number of
  class $\cF\cM$ and the same odd integer order of criticality.
  Suppose that the uniform distance
  $$
  \dist_{C^0}(\cR^nf,\cR^ng)\to 0\text{ geometrically fast}.
  $$
  Then there exists a $C^{1+\alpha}$ diffeomorphism of the circle
  conjugating $f$ and $g$.}

\medskip
\noindent
In combination with the main theorem of this paper, this implies:

\medskip
\noindent
{\bf Corollary.} {\sl Any two analytic critical circle mappings $f$ and $g$ with the same 
order of the critical point and the same irrational
rotation number $\rho\in\cF\cM$ are $C^{1+\alpha}$-smoothly conjugate.
}

\medskip
\noindent
de~Faria and de~Melo also provided examples of smooth
maps of unbounded type, which have the same rotation numbers and yet
are not $C^{1+\alpha}$ conjugate on the whole circle.
The first natural question to ask, given the main result of this paper is:

\medskip
\noindent
{\bf Question.} {\sl Does the conjugacy between two analytic critical circle maps 
with an irrational rotation number, which maps the critical point to the
critical point, have the smoothness $C^{1+\alpha}$ on the whole circle?}

\medskip
\noindent
It is also natural to ask:

\medskip
\noindent
{\bf Question.}  {\sl Are any two smooth critical circle maps
of bounded type  $C^{1+\alpha}$ conjugate at the critical point?}

\medskip
\noindent
 The counterexamples constructed by de~Faria and de~Melo make an essential 
use of the almost parabolic dynamics due to the unbounded type.

Further, there was a recent announcement by K.~Khanin and A.~Teplinskiy
of the following theorem:

\medskip

\noindent
{\bf Theorem \cite{Khanin}} {\sl Let $f$ and $g$ be two smooth critical circle mappings
with the
same irrational rotation number and the same odd integer order of criticality. Suppose
that the $C^{2+\eps}$-distance
$$\dist_{C^{2+\eps}}(\cR^nf,\cR^ng)\to 0\text{ geometrically fast}.$$
Then there exists a $C^1$-diffeomorphism of the circle conjugating $f$ and $g$.}

\medskip
\noindent
In combination with the renormalization convergence result of \cite{Ya4} (which we
proved by a different method in this paper), this implies:

\medskip
\noindent
{\bf Corollary.} {\sl Any two analytic critical circle mappings $f$ and $g$ with the same 
order of the critical point and the same irrational
rotation number are $C^1$-smoothly conjugate.
}

\medskip
\noindent
It is not clear if this result can be generalized to the smooth case. In particular,

\medskip
\noindent
{\bf Question.} {\sl Let $f$ and $g$ be two $C^3$-smooth critical circle maps with the same irrational
rotation number and the same order of the critical point.
 Is it always true that $$\dist(\cR^nf,\cR^ng)\to 0\text{ geometrically fast}?$$}

\medskip

Renormalization theory for analytic critical circle maps with a cubic critical point is now complete
\cite{Ya3,Ya4}, and all the results generalize to the case of a 
map with an arbitrary odd order $\gamma>1$ of the critical point.
However, numerical experiments show that the same kind of renormalization horseshoe
exists also for non odd-integer $\gamma>1$. The obvious issues with
analyticity preclude from directly extending the existing theory to
those cases. This remains therefore an important and challenging
problem.

\end{document}